\documentclass[12pt]{elsarticle}
\usepackage[letterpaper,top=2cm,bottom=2cm,left=3cm,right=3cm,marginparwidth=1.75cm]{geometry}
\usepackage{amsmath}
\usepackage{amsthm}
\usepackage{amsfonts}
\usepackage{graphicx}
\usepackage{longtable}
\usepackage{booktabs} 
\usepackage{algorithm}
\usepackage{algpseudocode}
\usepackage{multirow}
\usepackage{tikz}
\usepackage{subfig}
\usepackage{overpic}
\usepackage{wrapfig}
\usepackage{multicol}
\usepackage{multirow}

\usepackage[colorlinks=true, allcolors=blue]{hyperref}
\usepackage{url}
\usepackage{hyperref}

\usepackage{caption}
\graphicspath{{figures/}}
\makeatletter
\def\ps@pprintTitle{%
 \let\@oddhead\@empty
 \let\@evenhead\@empty
 \def\@oddfoot{Preprint submitted \it\hfill\today}%
 \let\@evenfoot\@oddfoot}
\makeatother

\begin{document}
\journal{}
\begin{frontmatter}
\title{Enhancing Arterial Blood Flow Simulations through Physics-Informed Neural Networks}
\author[1]{Shivam Bhargava}
\ead{shivam19@iisertvm.ac.in}
\author[1,2]{Nagaiah Chamakuri}
\ead{nagaiah.chamakuri@iisertvm.ac.in}
\address[1]{School of Mathematics, IISER Thiruvananthapuram, 695551 Kerala, India.}
\address[2]{Center for High-Performance Computing, IISER Thiruvananthapuram, 695551 Kerala, India}

\begin{abstract}
This study introduces a computational approach leveraging Physics-Informed Neural Networks (PINNs) for the efficient computation of arterial blood flows, particularly focusing on solving the incompressible Navier-Stokes equations by using the domain decomposition technique. Unlike conventional computational fluid dynamics methods, PINNs offer advantages by eliminating the need for discretized meshes and enabling the direct solution of partial differential equations (PDEs). In this paper, we propose the weighted Extended Physics-Informed Neural Networks (WXPINNs) and weighted Conservative Physics-Informed Neural Networks (WCPINNs), tailored for detailed hemodynamic simulations based on generalized space-time domain decomposition techniques. The inclusion of multiple neural networks enhances the representation capacity of the weighted PINN methods. Furthermore, the weighted PINNs can be efficiently trained in parallel computing frameworks by employing separate neural networks for each sub-domain. We show that PINNs simulation results circumvent backflow instabilities, underscoring a notable advantage of employing PINNs over traditional numerical methods to solve such complex blood flow models. They naturally address such challenges within their formulations. The presented numerical results demonstrate that the proposed weighted PINNs outperform traditional PINNs settings, where sub-PINNs are applied to each subdomain separately. This study contributes to the integration of deep learning methodologies with fluid mechanics, paving the way for accurate and efficient high-fidelity simulations in biomedical applications, particularly in modeling arterial blood flow.
\end{abstract}

\begin{keyword}
Physics-Informed Neural Networks, Hemodynamic Simulations, Incompressible Navier-Stokes, Domain decomposition, weighted XPINNs and weighted CPINNs. 
\end{keyword}

\end{frontmatter}
\section{Introduction}
The study of arterial blood flow occupies a pivotal junction in cardiovascular research, bridging advanced theoretical fluid dynamics with pivotal clinical applications. This interdisciplinary research is fundamental for deciphering the complex dynamics governing blood movement through arterial pathways, contributing significantly to our comprehension of the cardiovascular system. The precise simulation of arterial blood flow is not only critical for advancing theoretical understanding but also instrumental in improving diagnostic and therapeutic strategies in cardiovascular medicine \cite{Goubergrits_JMRI15, Schubert_SR20, Regitz_EHJ15}.

Central to the simulation of arterial blood flow is the incompressible Navier-Stokes equations, which present formidable challenges due to their nonlinear nature and the complexity of arterial geometries. In recent decades, researchers have extensively utilized various computational techniques to simulate and analyze blood flow, aiming to comprehend the correlation between vascular diseases and hemodynamics \cite{Purkayastha_JCBFM14}. Tokuda et al.~\cite{Tokuda_EJCTS08} employed the finite element method to numerically simulate blood flow, with a focus on understanding stroke mechanisms during cardio-pulmonary bypass. However, such numerical simulations of hemodynamics pose significant computational challenges in terms of time and memory requirements. To address this issue, numerous approaches have been explored to expedite the resolution of these problems \cite{Gerbeau_M2AN03, Crosetti_CF11, Wang_JSC18, Sarah_IJNMBE23, Badia_SISC08}. Such numerical methods for solving these equations fall short for numerous medical applications, like swiftly calculating hemodynamics for emergency cases and providing real-time surgical guidance. Moreover, in clinical practice, employing computational fluid dynamics necessitates repetitive simulations for varying patients, adding to the physician's workload. Given these challenges, there's a growing need to seek non-invasive, precise, cost-effective, and computationally efficient methods for capturing cardiovascular flow dynamics. Data-driven deep learning algorithms offer a promising avenue to address these issues. This necessitates the exploration of novel computational methodologies that can offer both accuracy and efficiency.

In this context, Physics-Informed Neural Networks (PINNs) emerge as a groundbreaking approach, melding the predictive power of machine learning with the foundational principles of physics \cite{raissi2017physicsI, raissi2017physicsII, RAISSI2019686, osti_1852843}. PINNs leverage the universal approximation capabilities of neural networks, constrained by the governing physical laws, to solve differential equations efficiently \cite{pinnreview, MAO2020112789, jin2021nsfnets, osti_1631426, RAO2020207}. This research innovates by applying PINNs to the domain of arterial blood flow, tackling the Navier-Stokes equations across varied geometries with unprecedented computational efficiency and model robustness. Many investigations have investigated employing deep learning for predicting hemodynamics. These studies encompass forecasting hemodynamic parameters or clinical metrics, local blood flow velocity, vessel cross-sectional area, and blood pressure \cite{taebi2022deep,arzani2021uncovering, jin2021nsfnets}. Until now, most of the literature concerning DL-driven flow modeling predominantly utilizes fully-connected neural network architectures \cite{arzani2021uncovering, jin2021nsfnets, eivazi2022physics}. Only a handful of publications explore the application of more intricate network structures, such as convolutional neural networks \cite{eichinger2020stationary, ma2021physics, gao2021phygeonet}. In computational fluid dynamics (CFD) simulations, encountering backflow at open boundaries poses a frequent challenge, often resulting in unphysical oscillations and instability issues, even at moderate Reynolds numbers, see \cite{Cristobal_JCP16, Moghadam_CM11} and referenced therein. Traditional numerical methods require stabilization terms to address these backflow instabilities.
Conversely, PINNs offer a different approach. They don't rely on stabilization terms to mitigate backflow instabilities and inherently tackle such challenges within their formulations.

The primary objective of this work is to demonstrate the efficacy of PINNs in simulating arterial blood flow, overcoming existing computational barriers while ensuring high fidelity in the simulations. For extensive simulations of Partial Differential Equations (PDEs), the high training cost associated with Physics-Informed Neural Networks (PINNs) makes tackling large-scale PDEs inherently expensive, impacting their efficiency compared to traditional numerical methods. Thus, there's a critical need to expedite the convergence of these models without compromising performance. Domain decomposition techniques, widely utilized in conventional numerical methods like finite difference, finite volume, and finite element methods, offer a promising solution. Here, the computational domain is subdivided into multiple subdomains, with interactions occurring solely at shared boundaries where continuity conditions are enforced. Within the realm of deep learning frameworks, the application of domain decomposition approaches in PINN is explored, notably in the conservative PINN (cPINN) method \cite{jagtap2020conservative} for conservation laws. Beyond cPINN, other strategies within the Scientific Machine Learning (SciML) domain include employing local neural networks on partitioned subdomains, as demonstrated by Li et al.\cite{Li_D3M_IEEE19}, who leveraged the variational principle. Similarly, Kharazmi et al. \cite{Kharazmi_CMAME21} proposed a variational PINN framework. Moreover, eXtended PINNs (XPINNs) \cite{jagtap2020extended} offer a comprehensive solution for solving various PDEs, exhibiting advantages akin to cPINN, such as employing separate neural networks in each subdomain, efficient hyper-parameter adjustment, facile parallelization, and substantial representation capacity. Further refining of the simulation accuracy and efficiency was studied in \cite{shukla2021parallel}.
In this study, we introduce weighted XPINN and CPINN methodologies to tackle hemodynamics model equations, aiming to enhance computational efficiency through domain decomposition strategies. By leveraging a finite number of subdomains, massively parallel computation becomes feasible, empowering effective handling of large-scale problems through domain decomposition. These weighted PINN approaches are more generalizations of the cPINNs \cite{jagtap2020conservative} and XPINNs\cite{jagtap2020extended}. Through a meticulous application of ADAM \cite{Adam} and L-BFGS \cite{LBFGS} optimization techniques and domain decomposition strategies, this research enhances the scalability and robustness of PINNs. These contributions mark significant advancements in integrating machine learning with fluid mechanics, offering a new paradigm for high-fidelity simulations in biomedical engineering and beyond.

The subsequent sections of this paper are structured as follows: Section 2 presents mathematical models concerning blood dynamics, encompassing the Navier-Stokes equations. Section 3 elaborates on computational methodologies, with particular emphasis on Physics-Informed Neural Networks (PINNs) and their diverse iterations. The application and refinement of these models in simulating blood flow dynamics are presented in Section 4. Following that, Section 5 unveils the outcomes and evaluations of the simulation performance. Lastly, Section 6 encapsulates the findings, addresses limitations, and delineates potential avenues for future research, highlighting the fusion of machine learning with fluid mechanics.

\section{Mathematical Formulation}
Blood displays viscoelastic properties and non-Newtonian behavior primarily at low shear rates due to the aggregation and alignment of erythrocytes, which enhance apparent viscosity and induce shear-thinning properties. Despite these complex behaviors, blood can be effectively modeled as a Newtonian fluid within arterial flow regimes, where shear rates are substantially higher. This model simplifies due to the uniform alignment of erythrocytes under high shear conditions, thereby simplifying the computational model without significantly sacrificing the accuracy of flow predictions in arterial geometries.

%\subsection{Formulation of Navier-Stokes Equations}
The dynamics of blood flow within the arterial system are governed by the incompressible Navier-Stokes equations, formulated within the computational domain \(\Omega \subset \mathbb{R}^2\) and expressed over the simulation interval \( (0, T] \). These equations capture the physiological nuances of arterial blood flow \cite{Quarteroni_CVS00}:
\begin{align}
    \rho \left( u_t + \mathbf{u} \cdot \nabla \mathbf{u} \right) + \nabla p - \mu \nabla^2 \mathbf{u} &= 0 \quad \text{in} \quad (0, T] \times \Omega, \\
    \nabla \cdot \mathbf{u} &= 0 \quad \text{in} \quad (0, T] \times \Omega,
\end{align}
where \(\mathbf{u}\) is the velocity field, \(p\) is the pressure field, \(\rho = 1060 \, \text{kg/m}^3\) is the blood density, and \(\mu = 3.5 \times 10^{-3} \, \text{Pa} \cdot \text{s}\) is the dynamic viscosity, assuming a Newtonian fluid approximation. \(T\) denotes the terminal simulation time.
%\textbf{write the density and viscosity coefficients????}

Incorporating the Cauchy stress tensor $\boldsymbol{\sigma}$, the equations can be alternatively formulated as:
\begin{equation}
    u_t + \mathbf{u} \cdot \nabla \mathbf{u} = \frac{1}{\rho} \nabla \cdot \boldsymbol{\sigma},
\end{equation}

with $\boldsymbol{\sigma} = -p\mathbf{I} + \mu (\nabla \mathbf{u} + (\nabla \mathbf{u})^T)$, and the pressure $p$ derived from:
\begin{equation}
    p = -\frac{1}{2} \text{tr}(\boldsymbol{\sigma}).
\end{equation}

For a two-dimensional flow, where velocity components are $u$ and $v$, the stress tensor $\boldsymbol{\sigma}$ is given by:
\begin{equation}
    \boldsymbol{\sigma} = 
    \begin{pmatrix}
    -p + 2\mu u_x & \mu (u_y + v_x) \\
    \mu (v_x + u_y) & -p + 2\mu v_y
    \end{pmatrix},
\end{equation}

{Initial and boundary conditions:}
To initiate the simulations, the velocity field throughout the computational domain is set to a quiescent state:
\[
\mathbf{u}(x, y, 0) = \begin{pmatrix} 0 \\ 0 \end{pmatrix}, \quad \forall (x, y) \in \Omega.
\]
This simplifies the analysis of fluid dynamics as the simulation progresses from a state of rest.

Boundary conditions are designed to mimic physiological conditions realistically and are adapted to the specific geometry and expected flow characteristics of each test case:

{Inlet conditions:}
At all inlet boundaries, the velocity is defined by a function that incorporates both the geometric specifics of the domain and the time-dependent aspects of physiological flow, often modulating according to a function that mimics cardiac pulsatility:
\[
\mathbf{u}(x_{\text{in}}, y, t) = U_{\text{max}} f(y, t) \begin{pmatrix} 1 \\ 0 \end{pmatrix},
\]
where \(x_{\text{in}}\) represents the inlet boundary coordinate, \(U_{\text{max}}\) is a scaling factor for maximum velocity, and \(f(y, t)\) is a profile function tailored to the specific domain geometry and flow conditions.

{Outlet and wall conditions:}
At the outlet, a condition is set to simulate a specific pressure condition depending on the study's need:
\[
p(x_{\text{out}}, y) = p_{\text{out}},
\]
where \(x_{\text{out}}\) is the outlet boundary coordinate, and \(p_{\text{out}}\) is typically set to zero or another value reflecting ambient or downstream pressure conditions.
%
%{Wall Conditions:}
No-slip boundary conditions are enforced along the walls of the domain to reflect the physical reality that the fluid velocity at a solid boundary is zero due to viscosity:
\[
\mathbf{u}(x, y_{\text{wall}}) = \begin{pmatrix} 0 \\ 0 \end{pmatrix},
\]
where \(y_{\text{wall}}\) denotes the coordinates of the wall boundary, applicable to both straight and curved segments in any domain shape.

These boundary conditions are crafted to ensure that the flow dynamics within the simulation accurately reflect the complex interactions expected in real-world fluid flow scenarios, especially within physiological contexts. They provide a robust framework that can adapt to both the specific details in the test cases developed in Section 4. 

\section{Computational Methodology}

Physics-Informed Neural Networks (PINNs) \cite{raissi2017physicsI,raissi2017physicsII,RAISSI2019686,osti_1852843} are designed to predict solutions \(\mathbf{u}(\mathbf{x}, t)\) for a set of partial differential equations (PDEs) across a specified domain \(\Omega\) and time interval \([0, T]\), where \(\mathbf{x}\) represents spatial coordinates and \(t\) total simulation time. The neural network, denoted by \(\mathcal{N}_{\boldsymbol{\theta}}(\mathbf{x}, t)\), with parameters \(\boldsymbol{\theta}\), is tasked with approximating the solution to the PDEs subject to initial and boundary conditions.

Given a differential operator \(\mathcal{D}\), a function \(f\) representing source terms, boundary conditions described by a function \(g\), and initial conditions defined by \(h\), the problem can be formally expressed as:
\begin{align*}
\mathcal{D}[\mathbf{u}](\mathbf{x}, t) &= f(\mathbf{x}, t), & \mathbf{x} &\in \Omega, \, t \in [0, T], \\
\mathcal{B}[\mathbf{u}](\mathbf{x}, t) &= g(\mathbf{x}, t), & \mathbf{x} &\in \partial\Omega, \, t \in [0, T], \\
\mathbf{u}(\mathbf{x}, 0) &= h(\mathbf{x}), & \mathbf{x} &\in \Omega.
\end{align*}

The core of the PINN methodology lies in the formulation of a physics-informed loss function, which quantifies the network's deviation not only from empirical data but also from the physical laws encapsulated by the governing differential equations. This loss function is constructed by evaluating the residuals associated with the governing PDEs (\(R_{\text{g}}\)), boundary conditions (\(R_{\text{bc}}\)), and initial conditions (\(R_{\text{ic}}\)) across the domain and its boundaries. These residuals are defined as follows:
\begin{align*}
R_{\text{g}}(\mathbf{x}, t; \boldsymbol{\theta}) &= \mathcal{D}[\mathcal{N}_{\boldsymbol{\theta}}](\mathbf{x}, t) - f(\mathbf{x}, t), \\
R_{\text{bc}}(\mathbf{x}, t; \boldsymbol{\theta}) &= \mathcal{B}[\mathcal{N}_{\boldsymbol{\theta}}](\mathbf{x}, t) - g(\mathbf{x}, t), \\
R_{\text{ic}}(\mathbf{x}; \boldsymbol{\theta}) &= \mathcal{N}_{\boldsymbol{\theta}}(\mathbf{x}, 0) - h(\mathbf{x}).
\end{align*}

The aggregated loss function for the weighted PINN, combining these residuals, is given by:
\begin{equation*}
\mathcal{L}_{\text{WPINN}}(\boldsymbol{\theta}) = \frac{1}{N_g}\sum_{j=1}^{N_g} \left\| R_{\text{g}}(\mathbf{x}_j, t_j; \boldsymbol{\theta}) \right\|^2, + \beta (\frac{1}{N_{bc}}\sum_{k=1}^{N_{bc}} \left\| R_{\text{bc}}(\mathbf{x}_k, t_k; \boldsymbol{\theta}) \right\|^2 + \frac{1}{N_{ic}}\sum_{l=1}^{N_{ic}} \left\| R_{\text{ic}}(\mathbf{x}_l; \boldsymbol{\theta}) \right\|^2)
\end{equation*}
where \(N_g\), \(N_{bc}\), and \(N_{ic}\) denote the number of collocation points for the governing equations, boundary conditions, and initial conditions, respectively. Here \(\beta>0\) is taken as a user-defined weighting coefficient that balances \( \mathcal{L}_{\text{g}} \) (governing equations loss) and \( \mathcal{L}_{\text{bc/ic}} \) (boundary and initial conditions loss) and accelerates convergence, enabling prioritization based on their physical relevance to enhance model accuracy and stability. 
Through the minimization of the above physics-informed loss, weighted PINNs effectively learn the dynamics of the system directly from the governing equations, offering a powerful tool for solving complex physical problems where traditional numerical methods may face limitations.

The techniques XPINNs\cite{jagtap2020extended} and CPINNs\cite{jagtap2020conservative} advance the capabilities of PINNs through domain decomposition. The main idea behind this methodology is partitioning the computational domain into multiple subdomains, $\Omega_i$, where each with a dedicated sub-network, $\mathcal{N}_{\theta_i}$, that facilitates localized solutions of complex PDEs.
%\subsubsection{Theoretical Framework}
For each subdomain $\Omega_i$, a specific sub-PINN $\mathcal{N}_{\theta_i}$ is trained to approximate the local solution $\mathbf{u}_i$. This process adheres to the governing differential equations, boundary conditions, and initial conditions particular to $\Omega_i$:
\begin{align}
\mathcal{D}_i[\mathbf{u}_i](\mathbf{x}, t; \boldsymbol{\theta}_i) &= f_i(\mathbf{x}, t), & \mathbf{x} \in \Omega_i, \, t \in [0, T], \\
\mathcal{B}_i[\mathbf{u}_i](\mathbf{x}, t; \boldsymbol{\theta}_i) &= g_i(\mathbf{x}, t), & \mathbf{x} \in \partial\Omega_i, \, t \in [0, T], \\
\mathbf{u}_i(\mathbf{x}, 0; \boldsymbol{\theta}_i) &= h_i(\mathbf{x}), & \mathbf{x} \in \Omega_i, \, t = 0.
\end{align}
These equations ensure the model captures the local physics accurately within each subdomain.

For XPINNs, ensuring smooth transitions between subdomains involves calculating an averaged solution at the interface points, $\Gamma_{ij}$, between adjacent subdomains $\Omega_i$ and $\Omega_j$. This averaging is represented as:

\begin{equation}
\mathbf{u}_{\text{avg},ij}(\mathbf{x}, t) = \frac{1}{2} \left( \mathcal{N}_{\boldsymbol{\theta}i}(\mathbf{x}, t) + \mathcal{N}_{\boldsymbol{\theta}_j}(\mathbf{x}, t) \right).
\end{equation}

The interface residual, $R_{\text{interface},i}$, measures the difference at these points between the sub-PINN's prediction for $\Omega_i$ and the averaged solution, ensuring interface continuity:

\begin{equation}
R_{\text{interface},i}(\mathbf{x}, t; \boldsymbol{\theta}_i, \boldsymbol{\theta}_j) = \mathcal{N}_{\boldsymbol{\theta}_i}(\mathbf{x}, t) - \mathbf{u}_{\text{avg},ij}(\mathbf{x}, t).
\end{equation}

Minimizing $R_{\text{interface},i}$ across interfaces achieves a cohesive global solution across the computational domain.

In weighted XPINNs (WXPINNs), the physics-informed loss for each subdomain includes weighted terms for governing equations, boundary conditions, initial conditions, and interfaces. The overall loss function is given by:
\begin{align}
\mathcal{L}_{\text{WXPINN}}(\boldsymbol{\theta}_i) = &\, \frac{1}{N_{g,i}}\sum_{j=1}^{N_{g,i}} \left\| R_{\text{g},i} \right\|^2 + \beta \left( \frac{1}{N_{bc,i}}\sum_{k=1}^{N_{bc,i}} \left\| R_{\text{bc}, i} \right\|^2 + \frac{1}{N_{ic,i}}\sum_{l=1}^{N_{ic,i}} \left\| R_{\text{ic},i} \right\|^2 \right) \nonumber \\
& + \gamma \sum_{j \in \mathcal{N}(i)} \frac{1}{N_{ij}} \sum_{\mathbf{x} \in \Gamma_{ij}} \left\| R_{\text{interface},i} \right\|^2.
\end{align}

where \( N_{ij}\) is the number of interface points between subdomains \(\Omega_i\) and \(\Omega_j\), and \(\mathcal{N}(i)\) represents the neighboring subdomains of \(\Omega_i\). Here, \(\beta>0\) serves as a user-defined weighting coefficient that tackles the balance of boundary and initial conditions loss, and \(\gamma>0\) serves as a user-defined weighting coefficient that targets the enhancement of interface continuity. \(\gamma\) controls the emphasis on the interface loss component, \(\mathcal{L}_{\text{interface}}\), which quantifies discrepancies and ensures smooth transitions at the inter-domain boundaries between adjacent sub-networks. Adjusting \(\gamma\) allows for fine-tuning the model's capacity to seamlessly integrate solutions across these interfaces, which is pivotal for achieving a coherent and accurate global solution. 

Through the minimization of the above physics-informed loss, weighted XPINNs solve complex PDEs across large and geometrically intricate domains by addressing them piecemeal yet in a manner that preserves the integrity and continuity of the overall solution.

\begin{algorithm}[H]
\caption{Training Process for Weighted XPINNs}
\begin{algorithmic}
\State Initialize neural network $\mathcal{N}_{\theta_i}$ for each subdomain $\Omega_i$
\For{each subdomain $\Omega_i$}
\State Generate collocation points for governing equations: $\{(\mathbf{x}_j, t_{j})\}_{j=1}^{N_{g,i}} \subset \Omega_i \times [0, T]$
\State Generate collocation points for boundary conditions: $\{(\mathbf{x}_k, t_{k})\}_{k=1}^{N_{bc,i}} \subset \partial\Omega_i \times [0, T]$
\State Generate collocation points for initial conditions: $\{(\mathbf{x}_l)\}_{l=1}^{N_{ic,i}} \subset \Omega_i$ at $t = 0$
\State Generate collocation points for interfaces with neighbors: $\{(\mathbf{x}_m, t_{m})\}_{m=1}^{N_{if,ij}} \subset \Gamma_{ij}$ for all $j \in \mathcal{N}(i)$
\EndFor
\While{not converged}
\For{each subdomain $\Omega_i$}
\State Compute residuals: $R_{\text{g},i}(\mathbf{x}, t; \boldsymbol{\theta}i)$, $R_{\text{bc},i}(\mathbf{x}, t; \boldsymbol{\theta}_i)$, and $R_{\text{ic},i}(\mathbf{x}; \boldsymbol{\theta}_i)$
\For{each neighboring subdomain $\Omega_j$}
\State Compute interface residuals $R_{\text{interface},ij}(\mathbf{x}, t; \boldsymbol{\theta}_i, \boldsymbol{\theta}_j)$
\EndFor
\State Compute physics-informed loss $\mathcal{L}_{\text{XPINN}}(\boldsymbol{\theta}_i)$:
\State $\mathcal{L}_{\text{WXPINN}}(\boldsymbol{\theta}_i) =  \mathcal{L}_{\text{g}}(\boldsymbol{\theta}_i) + \beta \mathcal{L}_{\text{bc/ic}}(\boldsymbol{\theta}_i) + \gamma \mathcal{L}_{\text{interface}}(\boldsymbol{\theta}_i)$
\State Update weights and biases $\boldsymbol{\theta}_i$ to minimize $\mathcal{L}_{\text{WXPINN}}(\boldsymbol{\theta}_i)$
\EndFor
\EndWhile
\State Assemble global solution $\mathbf{u}(\mathbf{x}, t)$ from local solutions $\mathbf{u}_i(\mathbf{x}, t; \boldsymbol{\theta}_i)$
\end{algorithmic}
\end{algorithm} 

In CPINNs, conservation across subdomain interfaces is ensured by aligning flux vectors, $\mathbf{F}_i$ and $\mathbf{F}_j$, from adjacent subdomains $\Omega_i$ and $\Omega_j$. The conservation condition at these interfaces, $\Gamma_{ij}$, is expressed as:

\begin{equation}
\left( \mathbf{F}_i(\mathbf{x}, t; \boldsymbol{\theta}_i) - \mathbf{F}_j(\mathbf{x}, t; \boldsymbol{\theta}_j) \right) \cdot \mathbf{n}_{ij} = 0 \,,
\end{equation}
where $\mathbf{n}_{ij}$ denotes the unit normal vector at $\Gamma_{ij}$, pointing from $\Omega_i$ to $\Omega_j$. This equation ensures the net flux of conserved quantities across the interface is zero, upholding the essential physical conservation laws.
The flux residual, quantifying the adherence to conservation laws at the interface $\Gamma_{ij}$, is defined as:
\begin{equation}
R_{\text{flux},ij}(\mathbf{x}, t; \boldsymbol{\theta}_i, \boldsymbol{\theta}_j) = \left( \mathbf{F}_i(\mathbf{x}, t; \boldsymbol{\theta}_i) - \mathbf{F}_j(\mathbf{x}, t; \boldsymbol{\theta}j) \right) \cdot \mathbf{n}_{ij}\,,
\end{equation}

In weighted CPINNs, the physics-informed loss for each subdomain incorporates residuals from the governing equations, boundary conditions, initial conditions, interfaces, and conservation of flux:
\begin{align}
\mathcal{L}_{\text{WCPINN}}(\boldsymbol{\theta}_i) = &\, \frac{1}{N_{g,i}}\sum_{j=1}^{N_{g,i}} \left\| R_{\text{g},i} \right\|^2 + \beta (\frac{1}{N_{bc,i}}\sum_{k=1}^{N_{bc,i}} \left\| R_{\text{bc},i} \right\|^2 \nonumber + \frac{1}{N_{ic,i}}\sum_{l=1}^{N_{ic,i}} \left\| R_{\text{ic},i} \right\|^2) \\ 
& + \gamma \sum_{j \in \mathcal{N}(i)} \frac{1}{N_{ij}} \sum_{\mathbf{x} \in \Gamma_{ij}} \left\| R_{\text{interface},i} \right\|^2 + \delta \sum_{j \in \mathcal{N}(i)} \frac{1}{N_{ij}} \sum_{\mathbf{x} \in \Gamma_{ij}} \left\| R_{\text{flux},i} \right\|^2
\end{align}

Here, \(\beta>0\) and \(\gamma>0\) serve as user-defined weighting coefficients that tackle the balance of boundary and initial conditions loss and tackle interface continuity, respectively, and \(\delta>0\) acts as a user-defined weighting coefficient that emphasizes the conservation of flux at the interfaces between subdomains. This coefficient is critical for ensuring that key conservation properties, such as mass and momentum in the case of the incompressible Navier-Stokes equations, are maintained across subdomain boundaries. By tuning \(\delta\), the model can more effectively handle fluid dynamics problems where the accurate representation of flow and pressure continuity is essential for achieving realistic and physically accurate simulations. Adjusting \(\delta\) allows the model to precisely manage how these conservation laws influence the overall solution, ensuring that the continuity of physical quantities is preserved.

Through the minimization of the above physics-informed loss, weighted CPINNs solve complex PDEs across large and geometrically intricate domains by addressing them piecemeal yet in a manner that preserves the integrity and continuity of the overall solution.

\begin{algorithm}[H]
\caption{Training Process for Weighted CPINNs}
\begin{algorithmic}
\State Initialize neural network $\mathcal{N}_{\theta_i}$ for each subdomain $\Omega_i$
\For{each subdomain $\Omega_i$}
   \State Generate collocation points for governing equations: $\{(\mathbf{x}_j, t_{j})\}_{j=1}^{N_{g,i}} \subset \Omega_i \times [0, T]$
   \State Generate collocation points for boundary conditions: $\{(\mathbf{x}_k, t_{k})\}_{k=1}^{N_{bc,i}} \subset \partial\Omega_i \times [0, T]$
   \State Generate collocation points for initial conditions: $\{(\mathbf{x}_l)\}_{l=1}^{N_{ic,i}} \subset \Omega_i$ at $t = 0$
   \State Generate collocation points for interfaces with neighbors: $\{(\mathbf{x}_m, t_{m})\}_{m=1}^{N_{if,ij}} \subset \Gamma_{ij}$ for all $j \in \mathcal{N}(i)$
\EndFor
\While{not converged}
\For{each subdomain $\Omega_i$}
   \State Compute residuals: $R_{\text{g},i}(\mathbf{x}, t; \boldsymbol{\theta}i)$, $R_{\text{bc},i}(\mathbf{x}, t; \boldsymbol{\theta}_i)$, and $R_{\text{ic},i}(\mathbf{x}; \boldsymbol{\theta}_i)$
   \For{each neighboring subdomain $\Omega_j$}
       \State Compute interface residuals $R_{\text{interface},ij}(\mathbf{x}, t; \boldsymbol{\theta}_i, \boldsymbol{\theta}_j)$
       \State Compute flux residuals $R_{\text{flux},ij}(\mathbf{x}, t; \boldsymbol{\theta}_i, \boldsymbol{\theta}_j)$
   \EndFor
   \State Compute physics-informed loss $\mathcal{L}_{\text{WCPINN}}(\boldsymbol{\theta}_i)$:
   \State $\mathcal{L}_{\text{WCPINN}}(\boldsymbol{\theta}_i) =  \mathcal{L}_{\text{g}}(\boldsymbol{\theta}_i) + \beta \mathcal{L}_{\text{bc/ic}}(\boldsymbol{\theta}_i) + \gamma \mathcal{L}_{\text{interface}}(\boldsymbol{\theta}_i) + \delta \mathcal{L}_{\text{flux}}(\boldsymbol{\theta}_i)$
   \State Update weights and biases $\boldsymbol{\theta}_i$ to minimize $\mathcal{L}_{\text{WCPINN}}(\boldsymbol{\theta}_i)$
\EndFor
\EndWhile
\State Assemble global solution $\mathbf{u}(\mathbf{x}, t)$ from local solutions $\mathbf{u}_i(\mathbf{x}, t; \boldsymbol{\theta}_i)$
\end{algorithmic}
\end{algorithm}

\section{Model Implementation}

To apply weighted Physics-Informed Neural Networks (WPINNs) to the incompressible Navier-Stokes equations, a neural network that maps the spatiotemporal variables $\{t, \mathbf{x}\}^T$ to the mixed-variable solution $\{\psi, p, \boldsymbol{\sigma}\}$ is constructed. The input $\mathbf{x}$ has 2 components $(x, y)$. The output $\boldsymbol{\sigma}$ has three components $\boldsymbol{\sigma}=({\sigma}_{11}, {\sigma}_{12},{\sigma}_{22})^T$.

\begin{figure}[H]
\centering
\includegraphics[width=0.7\textwidth]{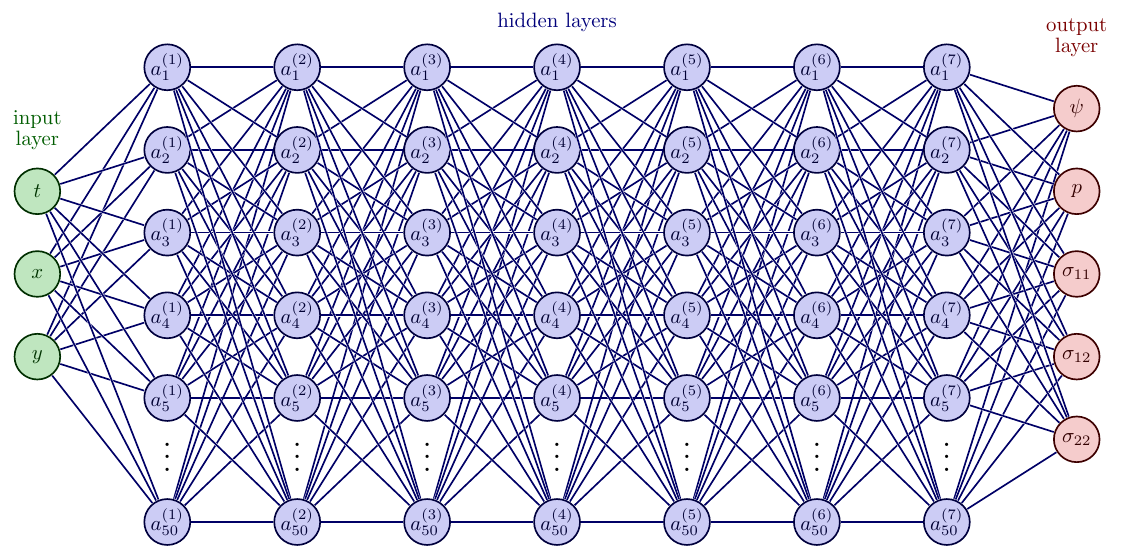} 
\caption{Representative Diagram of the PINN to solve Navier-Stokes Equations}
\end{figure}

The stream function $\psi$ is used rather than velocity $\mathbf{u}$ directly to ensure that the divergence-free condition of the flow is maintained. The velocity $\mathbf{u} = (u', v')$ where $u'$ and $v'$ are predicted velocities can be computed as:

\begin{equation}
    u' = \psi_y, \quad v' = -\psi_x
\end{equation}

Thus the residuals obtained from $u$ and $v$ can be computed as:

\begin{align}
\begin{pmatrix} R_{u} \\ R_{v} \end{pmatrix} = 
\begin{pmatrix} 
\rho u'_t + \rho (u' u'_x + v' u'_y) - \left(-p_x + 2\mu u'_{xx}\right) - \left(\mu (u'_y + v'_x)_y\right) \\
\rho v'_t + \rho (u' v'_x + v' v'_y) - \left(\mu (v'_x + u'_y)_x\right) - \left(-p_y + 2\mu v'_{yy}\right)
\end{pmatrix}
\end{align}

The residual obtained from $p$ for predicted pressure $p'$ can be computed as:

\begin{equation}
    R_{p} = p' + \frac{1}{2}({\sigma}_{11} + {\sigma}_{22})
\end{equation}

The residual obtained from ${\sigma}_{11}$, ${\sigma}_{12}$, ${\sigma}_{22}$ can be computed as:

\begin{equation}
     \begin{pmatrix} R_{{\sigma}_{11}} \\ R_{{\sigma}_{12}} \\ R_{{\sigma}_{22}} \end{pmatrix} = \begin{pmatrix} (-p' + 2\mu u'_x) - {\sigma}_{11} \\ (\mu u'_y + \mu v'_x) - {\sigma}_{12} \\ (-p' + 2\mu v'_y) - {\sigma}_{22} \end{pmatrix}
\end{equation}

The residual from the governing equations is given by:
\begin{equation}
    \|{R}_{\text{g}}(x, y, t)\|^2 = \|{R}_{u}\|^2 + \|{R}_{v}\|^2 + \|R_{p}\|^2 + \|{R}_{{\sigma}_{11}}\|^2 + \|{R}_{{\sigma}_{12}}\|^2 + \|{R}_{{\sigma}_{22}}\|^2
\end{equation}

The residuals from the boundary and initial conditions are given by:
\begin{align}
    \|{R}_{\text{bc}}(x, y, t)\|^2 &= \|{u}(x, y, t) - {u'}(x, y, t)\|^2 + \|{v}(x, y, t) - {v'}(x, y, t)\|^2 \nonumber \\ 
    &+ \|{p}(x, y, t) - {p'}(x, y, t)\|^2 \\
    \|{R}_{\text{ic}}(x, y, 0)\|^2 &= \|{u}(x, y, 0) - {u'}(x, y, 0)\|^2 + \|{v}(x, y, 0) - {v'}(x, y, 0)\|^2 \nonumber \\ 
    &+ \|{p}(x, y, 0) - {p'}(x, y, 0)\|^2
\end{align}

The physics-informed loss and initial/boundary condition loss can now be given by:
\begin{align}
   \mathcal{L}_{\text{g}} &= \frac{1}{N_g}\sum_{i=1}^{N_g} \| R_{g}(x_i, y_i, t_i) \|^2 \\
   \mathcal{L}_{\text{bc/ic}} &= \frac{1}{N_{bc}}\sum_{i=1}^{N_{bc}} \| R_{bc}(x_i, y_i, t_i) \|^2 + \frac{1}{N_{ic}}\sum_{i=1}^{N_{ic}} \| R_{ic}(x_i, y_i, 0) \|^2 
\end{align}
where \(N_{(\cdot)}\) represents the number of collocation points.

The loss function for the Physics-Informed Neural Network (PINN) model is formulated as follows:
\begin{equation}
    \mathcal{L}_{\text{WPINN}} = \mathcal{L}_{\text{g}} + \beta \mathcal{L}_{\text{bc/ic}},
\end{equation}
where $\beta$ serves as a weighting factor to balance the components of the loss function.

To apply the weighted XPINN to the incompressible Navier-Stokes equations, the domain $\Omega$ is decomposed into $M$ non-overlapping subdomains $\{\Omega_i\}_{i=1}^M$, each with its own sub-PINN $\mathcal{N}_i$. The interface between adjacent subdomains $\Omega_i$ and $\Omega_{j}$ is denoted as $\Gamma_{i,j}$ where $j \in \{i-1, i+1\}$. For each subdomain $\Omega_i$, interface conditions are enforced at the shared boundaries with adjacent subdomains. The predicted velocity fields $u_{i}'$ and $v_{i}'$ and the pressure field $p_i'$ from the sub-PINN for $\Omega_i$, and $u_{j}'$, $v_{j}'$, and $p_j'$ from $\Omega_j$ are used to define the interface residuals.

The interface residual for velocity and pressure for subdomain $\Omega_i$ at the interface $\Gamma_{i,j}$ is given by:

\begin{align}
R_{\text{interface},u}(x, y, t) &= u_{i}'(x, y, t) - \frac{u_{i}'(x, y, t) + u_{j}'(x, y, t)}{2}, \\
R_{\text{interface},v}(x, y, t) &= v_{i}'(x, y, t) - \frac{v_{i}'(x, y, t) + v_{j}'(x, y, t)}{2}, \\
R_{\text{interface},p}(x, y, t) &= p_i'(x, y, t) - \frac{p_i'(x, y, t) + p_{j}'(x, y, t)}{2},
\end{align}

The interface loss $\mathcal{L}_{\text{interface},ij}$ for subdomain $\Omega_i$ for the interface $\Gamma_{i,j}$ is calculated based on the interface residuals:

\begin{equation}
\scalebox{0.9}{$\displaystyle \mathcal{L}_{\text{interface},ij} = \frac{1}{N_{ij}}\sum_{k=1}^{N_{ij}} \left( \left\| R_{\text{interface},u}(x_k, y_k, t_k) \right\|^2 + \left\| R_{\text{interface},v}(x_k, y_k, t_k) \right\|^2 + \left\| R_{\text{interface},p}(x_k, y_k, t_k) \right\|^2 \right)$}
\end{equation}

Thus the interface loss $\mathcal{L}_{\text{interface},i}$ for subdomain $\Omega_i$ becomes:
\begin{equation}
    \mathcal{L}_{\text{interface},i} = \mathcal{L}_{\text{interface},i(i+1)} + \mathcal{L}_{\text{interface},i(i-1)} \,.
\end{equation}

The total loss function for each sub-PINN $\mathcal{N}_i$ in subdomain $\Omega_i$ is a combination of the physics-informed loss, boundary, and initial condition loss, and the interface loss:

\begin{equation}
\mathcal{L}_{\text{WXPINN},i} =  \mathcal{L}_{\text{g},i} + \beta \mathcal{L}_{\text{bc/ic},i} + \gamma \mathcal{L}_{\text{interface},i} 
\end{equation}
where $\gamma$ serves as a weighting factor to balance the components of the loss function.

%\subsection{Weighted CPINN for Navier-Stokes Equations}
To apply the weighted CPINN to the incompressible Navier-Stokes equations, the domain $\Omega$ is decomposed into $M$ non-overlapping subdomains $\{\Omega_i\}_{i=1}^M$, each with its own sub-PINN $\mathcal{N}_i$. The interface between adjacent subdomains $\Omega_i$ and $\Omega_{j}$ is denoted as $\Gamma_{i,j}$ where $j \in \{i-1, i+1\}$. For each subdomain $\Omega_i$, interface conditions are enforced at the shared boundaries with adjacent subdomains. The predicted velocity fields $u_{i}'$ and $v_{i}'$ and the pressure field $p_i'$ from the sub-PINN for $\Omega_i$, and $u_{j}'$, $v_{j}'$, and $p_j'$ from $\Omega_j$ are used to define the flux residuals for conservation of momentum and conservation of energy.

The flux residuals for momentum and energy for subdomain $\Omega_i$ at the interface $\Gamma_{i,j}$ are given by:

\begin{align}
R_{\text{flux},m}(x, y, t) &= \rho \left( u_{i}'(x, y, t) - u_{j}'(x, y, t) \right), \\
R_{\text{flux},\mu}(x, y, t) &= \left(\rho  u_{i}'(x, y, t)^2 + p_{i}'(x, y, t) \right) - \left(\rho u_{j}'(x, y, t)^2 + p_{j}'(x, y, t) \right),
\end{align}

The flux loss $\mathcal{L}_{\text{f},ij}$ for subdomain $\Omega_i$ for the interface $\Gamma_{i,j}$ is calculated based on the flux residuals:

\begin{equation}
 \mathcal{L}_{\text{flux},ij} = \frac{1}{N_{ij}}\sum_{k=1}^{N_{ij}} \left( \left\| R_{\text{flux},m}(x_k, y_k, t_k) \right\|^2 + \left\| R_{\text{flux},\mu}(x, y, t)(x_k, y_k, t_k) \right\|^2 \right)
\end{equation}

Thus the flux loss $\mathcal{L}_{\text{flux},i}$ for subdomain $\Omega_i$ becomes:
\begin{equation}
    \mathcal{L}_{\text{flux},i} = \mathcal{L}_{\text{flux},i(i+1)} + \mathcal{L}_{\text{flux},i(i-1)}
\end{equation}

The total loss function for each sub-PINN $\mathcal{N}_i$ in subdomain $\Omega_i$ is a combination of the physics-informed loss, boundary, and initial condition loss, interface loss, and the flux loss:

\begin{equation}
\mathcal{L}_{\text{WCPINN},i} =  \mathcal{L}_{\text{g},i} + \beta \mathcal{L}_{\text{bc/ic},i} + \gamma \mathcal{L}_{\text{interface},i} + \delta \mathcal{L}_{\text{flux}, i}
\end{equation}

where $\delta$ serves as a weighting factor to balance the components of the loss function.

%\subsection{Test Cases} 
This study encompasses two primary test cases designed to demonstrate the efficacy of the PINN framework in simulating non-Newtonian fluid flow dynamics: a rectangular domain and a semi-circular domain.

Rectangular Domain:
Numerical Simulation of flow in a rectangle using PINNs is presented. The computational domain is a rectangle in two dimensions: a length of 1.1 cm and a height of 0.41 cm. The viscosity \(\nu\) is defined as 0.01 gs\(^{-1}\)cm\(^{-1}\). The total simulation time, \(T\), is 0.5 seconds with a time step \(\Delta t\) of 0.01 seconds.

%Initial Conditions:
The initial velocity field at \( t = 0 \) seconds is given by:
\[
\mathbf{u}(x, y, 0) = 
\begin{pmatrix}
    u(x, y, 0) \\
    v(x, y, 0)
\end{pmatrix}
=
\begin{pmatrix}
    0 \\
    0
\end{pmatrix}
\]
\noindent
for all points \((x, y)\) within the domain.

%Boundary Conditions:
The boundary conditions are given by:
\begin{itemize}
    \item Inlet (x = 0):
     A parabolic profile gives the velocity vector at the inlet boundary:
    \[
    \mathbf{u}(0,y) = 
    \begin{pmatrix}
        {u}(0, y) \\
        {v}(0, y)
    \end{pmatrix}
    =
   \scalebox{0.9}{$\displaystyle \begin{pmatrix}
        4{U}_{\text max}\dfrac{y(H - y)}{H^2}(\sin(\dfrac{\pi t}{T} + \dfrac{3\pi}{2}) + 1) \\
        0
    \end{pmatrix} $} 
    \]

    where \(H\) is the height of the rectangle and \({U}_{\text max}=0.5\).

    \item Outlet (x = L):
    At the outlet, a pressure condition:
    \[ p(L, y) = 0 \]

    Where \(L\) is the length of the rectangle.
   
    \item Solid Walls (y = 0 and y = H):
    No-slip conditions at the solid walls:
    \[
    \begin{aligned}
        {u}(x, 0) &= 0, & {v}(x, 0) &= 0, \\
        {u}(x, H) &= 0, & {v}(x, H) &= 0,
    \end{aligned}
    \]
\end{itemize}
\noindent

Our implementation of the weighted PINNs is based on using the TensorFlow framework. The neural network architecture consists of an input layer with three neurons, an output layer with five neurons, and seven hidden layers containing 50 neurons. The tanh function is utilized as the activation function for each hidden layer. A total of 3321 collocation points are generated for network training, which includes 244 boundary points (\(N_b\)) and 81 points on the inlet/outlet boundary (\(N_{in/out}\)), using Latin hypercube sampling (LHS). The balancing coefficient \(\beta\) is taken as 1, 2, 5, and 10.

For optimization, the network employs both Adam and L-BFGS optimizers. The network undergoes an initial training phase with the Adam optimizer for 5000 iterations, followed by a subsequent stage with the L-BFGS optimizer until convergence. The Hager-Zhang line search algorithm is utilized with a maximum of 50 iterations to determine the optimal step length \(\alpha_k\).

For prediction, 64561 collocation points are used, which includes 1124 boundary points (\(N_b\)) and 161 points on the inlet/outlet boundary (\(N_{in/out}\)).

For the WXPINN model, the computational domain was divided into \(M= \) 2, 3, and 4 non-overlapping subdomains. Each subdomain was modeled using a separate sub-PINN, whose network architecture mirrored that of the global model, and the interface conditions were enforced as described in the previous sections. The same number of collocation points as the global model were taken and equally distributed across the subdomains. The balancing coefficient \(\gamma\) is taken as 1, 2, 5, and 10 while keeping \(\beta\) fixed at 1. 

For the WCPINN model, the computational domain was divided into \(M= \) 2, 3, and 4 non-overlapping subdomains. Each subdomain was modeled using a separate sub-PINN, whose network architecture mirrored that of the global model, and the interface conditions were enforced as described in the previous sections. The same number of collocation points as the global model were taken and equally distributed across the subdomains. The balancing coefficient \(\delta\) is taken as 1, 2, 5, and 10 while keeping \(\beta\) fixed at 1 and \(\gamma\) taken as 1 and 5. 

Semi-Circular Domain: 
Numerical Simulation of flow in a semi-circular domain using PINNs is presented. The computational domain is a semi-circular pipe with smooth disturbances in two dimensions, possessing a cross-sectional radius \( a = 1.6 \) cm and curvature radius \( R = 2.9 \) cm. The kinematic viscosity \( \nu \) is set to \( 0.4 \) gs\(^{-1}\) cm\(^{-1}\). The total simulation time, \( T \), is 6 seconds with a time step \( \Delta t \) of 0.01 seconds.

%Initial Conditions:
The initial velocity field at \( t = 0 \) seconds is given by:
\[
\mathbf{u}(x, y, 0) = 
\begin{pmatrix}
    u(x, y, 0) \\
    v(x, y, 0)
\end{pmatrix}
=
\begin{pmatrix}
    0 \\
    0
\end{pmatrix}
\]
for all points \((x, y)\) within the domain.

%Boundary Conditions:
The boundary conditions are given by:
\begin{itemize}
    \item Inlet (at the straight segment):
    A parabolic profile gives the velocity vector at the inlet boundary:
    \[
    \mathbf{u}_{\text{inlet}} = 
    \begin{pmatrix}
        {u}_{\text{inlet}} \\
        {v}_{\text{inlet}}
    \end{pmatrix}
    =
    \scalebox{0.9}{$\displaystyle \begin{pmatrix}
        4{U}_{\text{max}}\dfrac{x(D - x)}{D^2}(\sin(\dfrac{\pi t}{T} + \dfrac{3\pi}{2}) + 1) \\
        0
    \end{pmatrix} $} 
    \]
    where \( U_{\text{max}} = 0.75\) is the maximum velocity at the inlet, and \( D \) is the cross-sectional diameter of the inlet segment.

    \item Outlet (at the straight segment):
    At the outlet, a pressure condition:
    \[ p(x_{\text{out}}, y) = 0 \]
    where \( x_{\text{out}} \) denotes the position at the outlet.

    \item Curved Boundary (semi-circular edge):
    No-slip conditions are applied at the curved boundary:
    \[
    \mathbf{u}_{\text{curved}} = 
    \begin{pmatrix}
        {u}_{\text{curved}} \\
        {v}_{\text{curved}}
    \end{pmatrix}
    =
    \begin{pmatrix}
        0 \\
        0
    \end{pmatrix}
    \]
\end{itemize}

The weighted PINN model is implemented in TensorFlow. The neural network architecture consists of an input layer with three neurons, an output layer with five neurons, and seven hidden layers containing 50 neurons. The tanh function is utilized as the activation function for each hidden layer. In total, 29760 collocation points are employed, with a batch size of 20000 points used in each training iteration. The balancing coefficient \(\beta\) is set to 1, 2, 5, and 10.

For optimization, the network employs both Adam and L-BFGS optimizers. The network undergoes an initial training phase with the Adam optimizer for the first few iterations, followed by a subsequent stage with the L-BFGS optimizer until convergence. The Hager-Zhang line search algorithm is utilized with a maximum of 50 iterations to determine the optimal step length \(\alpha_k\). For prediction, the exact total of 29760 collocation points is utilized. 
The numerical results showcased in this study were obtained using a Linux cluster located at IISER Thiruvananthapuram. The cluster comprises 88 nodes, with each node housing 28 compute cores running at 2.60 GHz and equipped with 128GB of RAM.

\section{Results}
%\subsection{PINN Model}
%\subsubsection{Rectangular Domain}
The weighted PINN model's capability to simulate fluid dynamics within a rectangular domain was evaluated through the visualization of velocity and pressure fields at various timestamps. These visualizations illustrated the dynamic flow patterns, identifying regions of varying velocity and pressure essential for understanding the fluid's behavior and the impact of forces on domain boundaries.

\begin{figure}[H]
    \centering
    \includegraphics[width=0.325\textwidth]{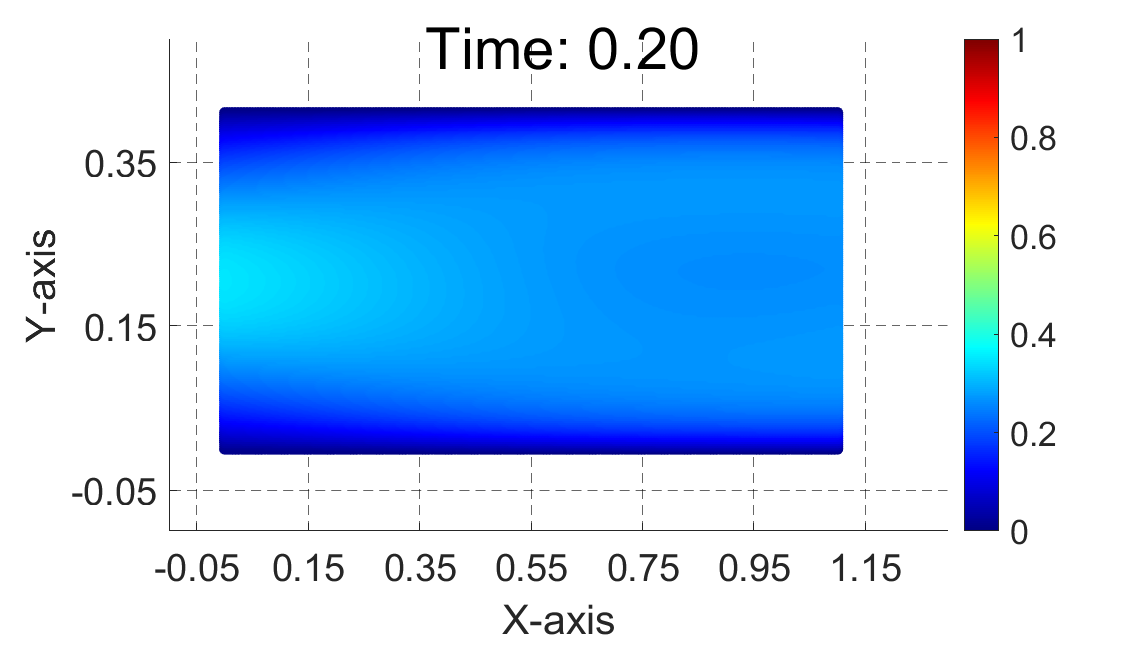}
    \includegraphics[width=0.325\textwidth]{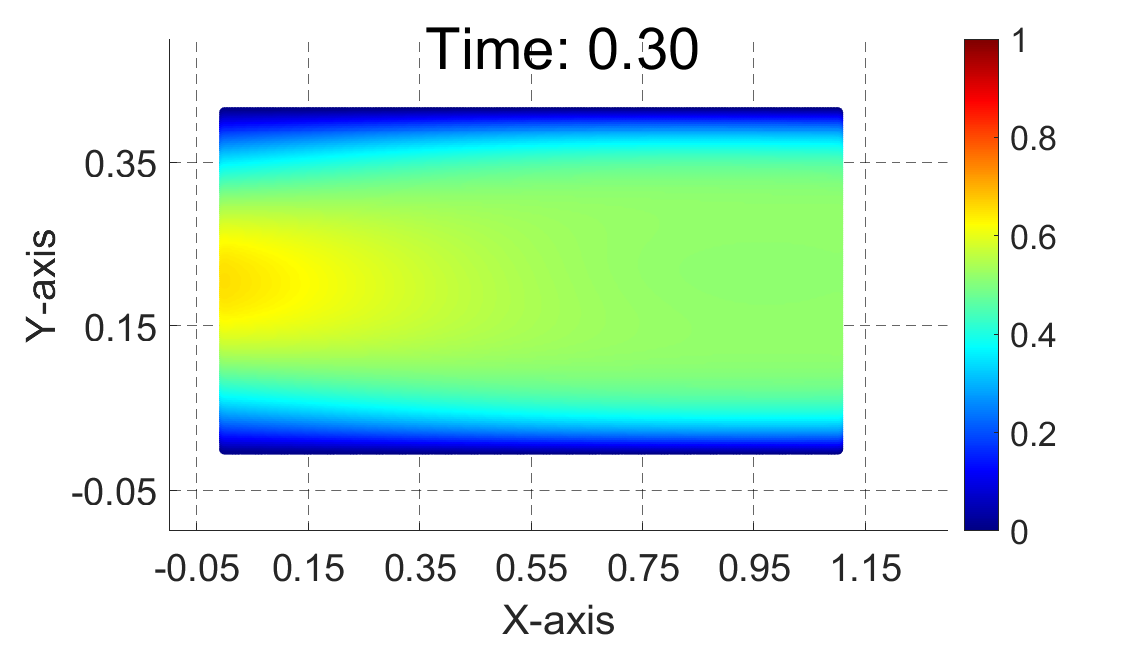}
    \includegraphics[width=0.325\textwidth]{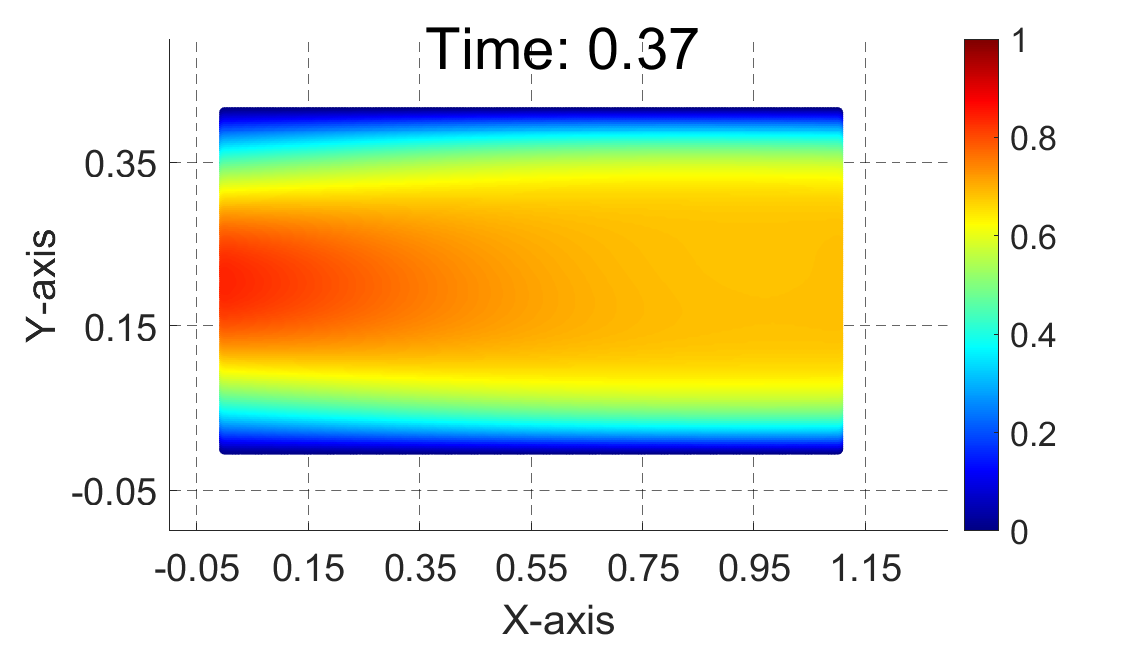}
    \includegraphics[width=0.325\textwidth]{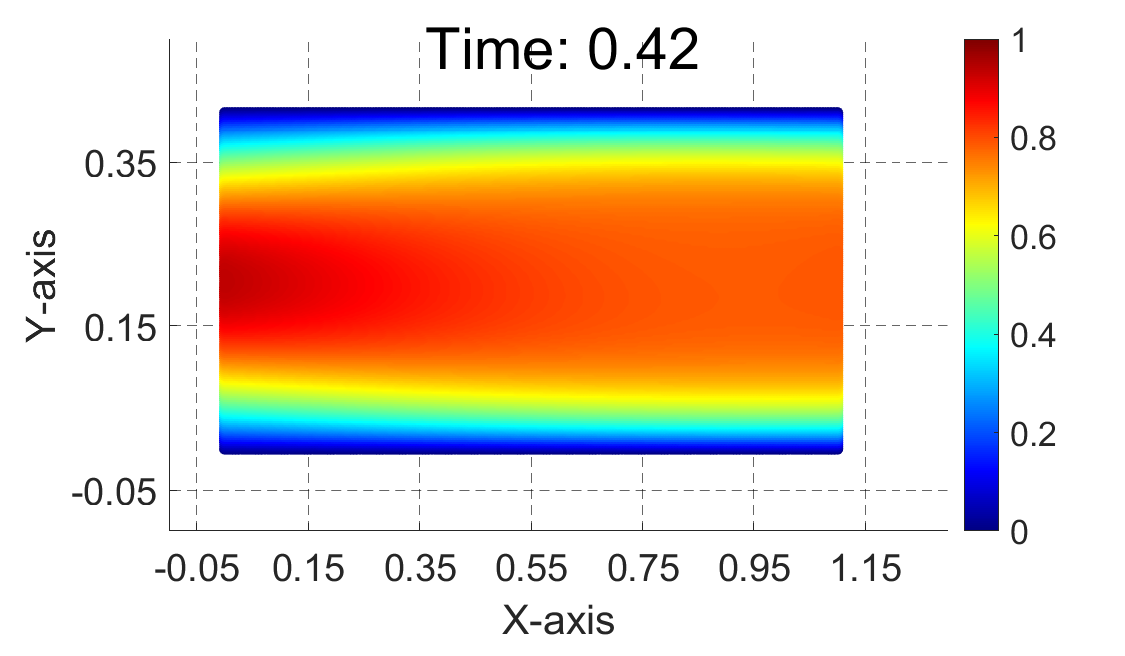}
    \includegraphics[width=0.325\textwidth]{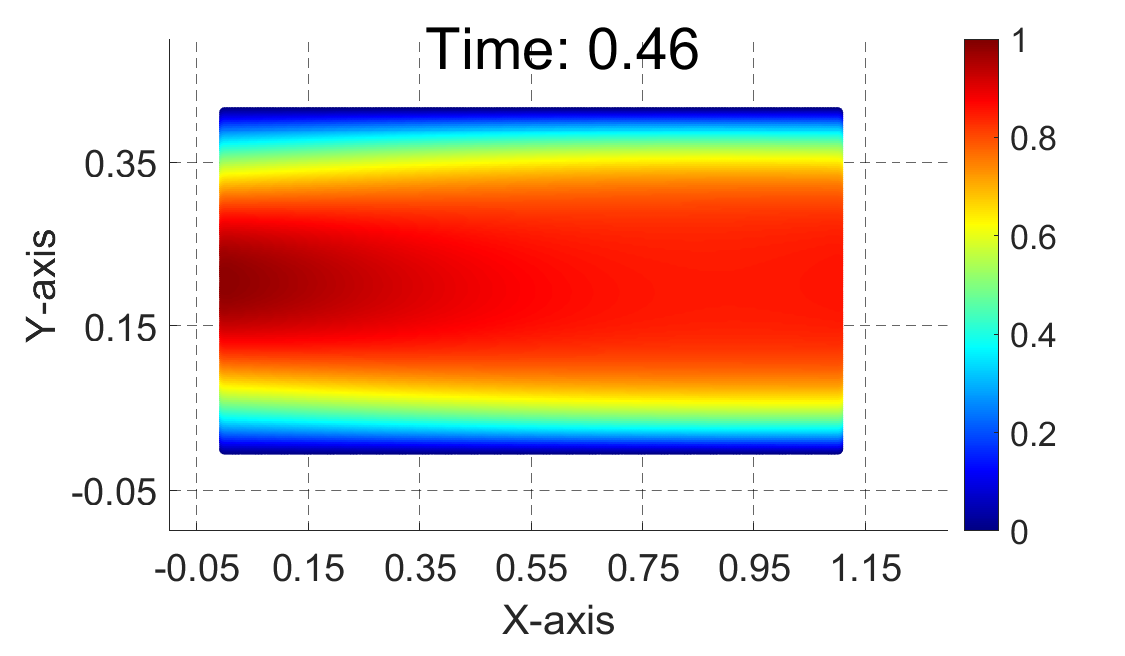}
    \includegraphics[width=0.325\textwidth]{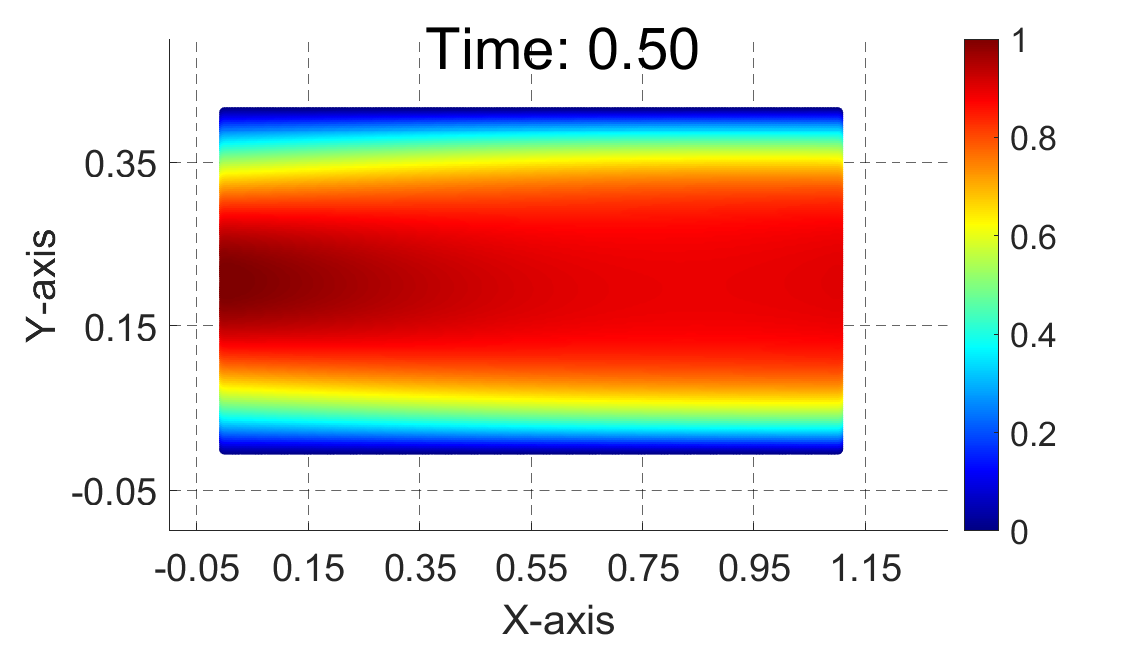}
    \caption{Predicted velocity field at various timestamps obtained from the weighted PINN model for the rectangular domain}
    \label{fig:rect_vel_field-1}
\end{figure}

%Velocity Field Analysis: 
The model effectively captured the evolution of the flow field over time, highlighting significant variations in velocity across the domain, which is depicted in Figure~\ref{fig:rect_vel_field-1} for the rectangular domain. The velocity profile is consistent, with the maximum velocity observed in the middle of the inlet decreasing towards the walls due to the no-slip boundary conditions. The velocity in the middle of the channel is observed to continuously increase in the simulation time period due to the parabolic inlet profile.

%Predicted Velocity Plots for rectangular domain:

\begin{figure}[H]
    \centering
    \includegraphics[width=0.325\textwidth]{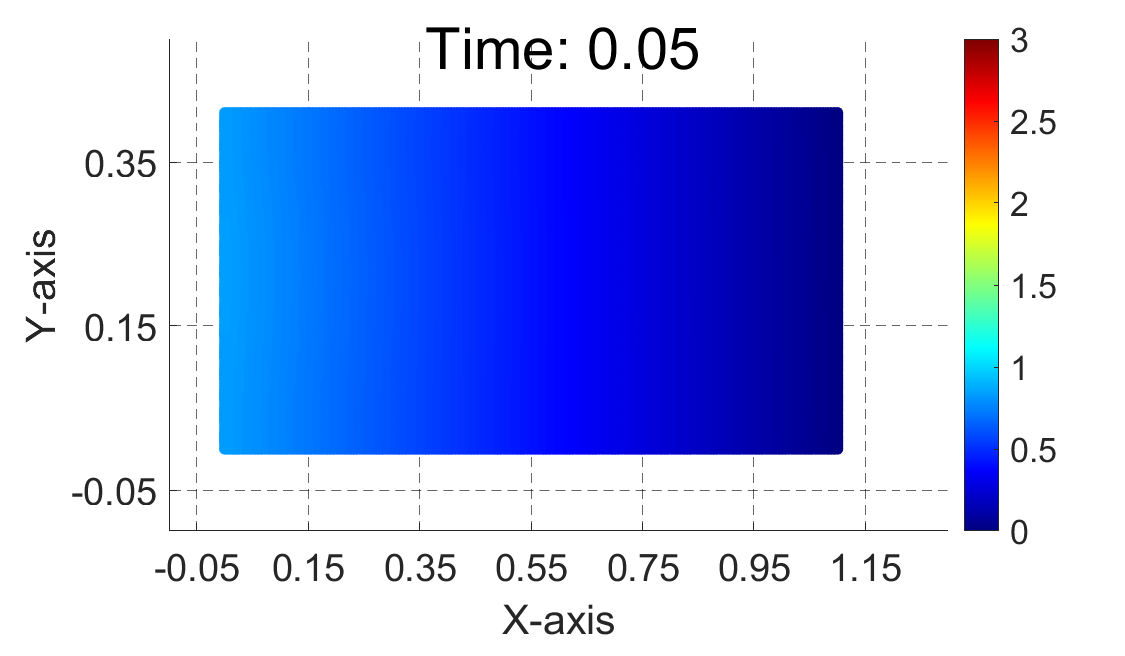}
    \includegraphics[width=0.325\textwidth]{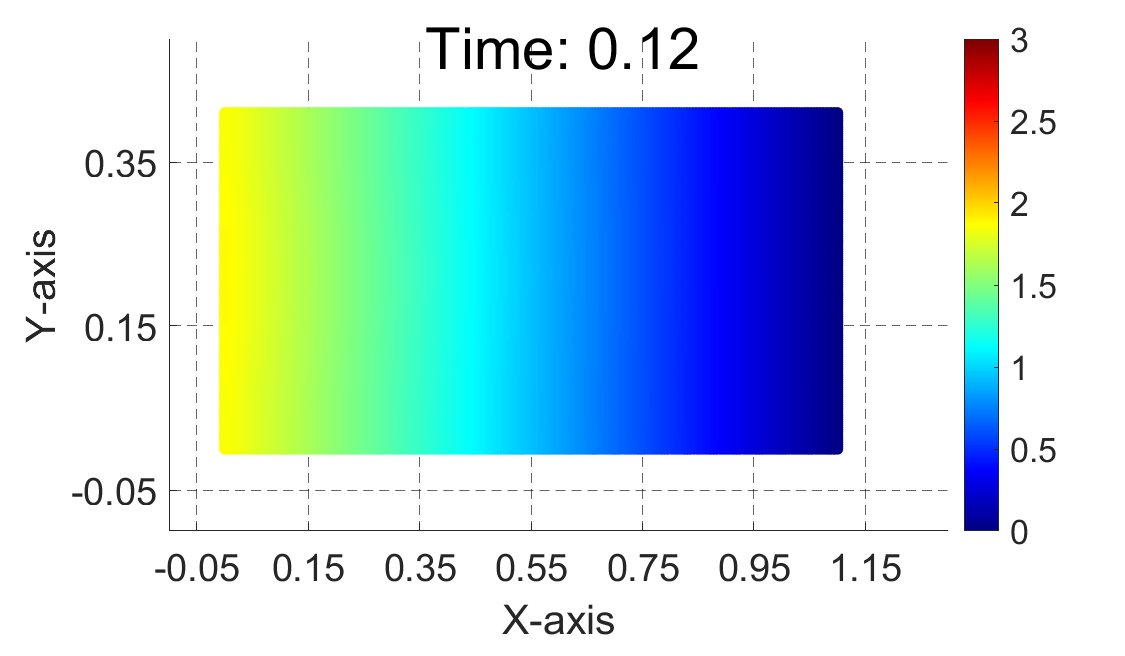}
    \includegraphics[width=0.325\textwidth]{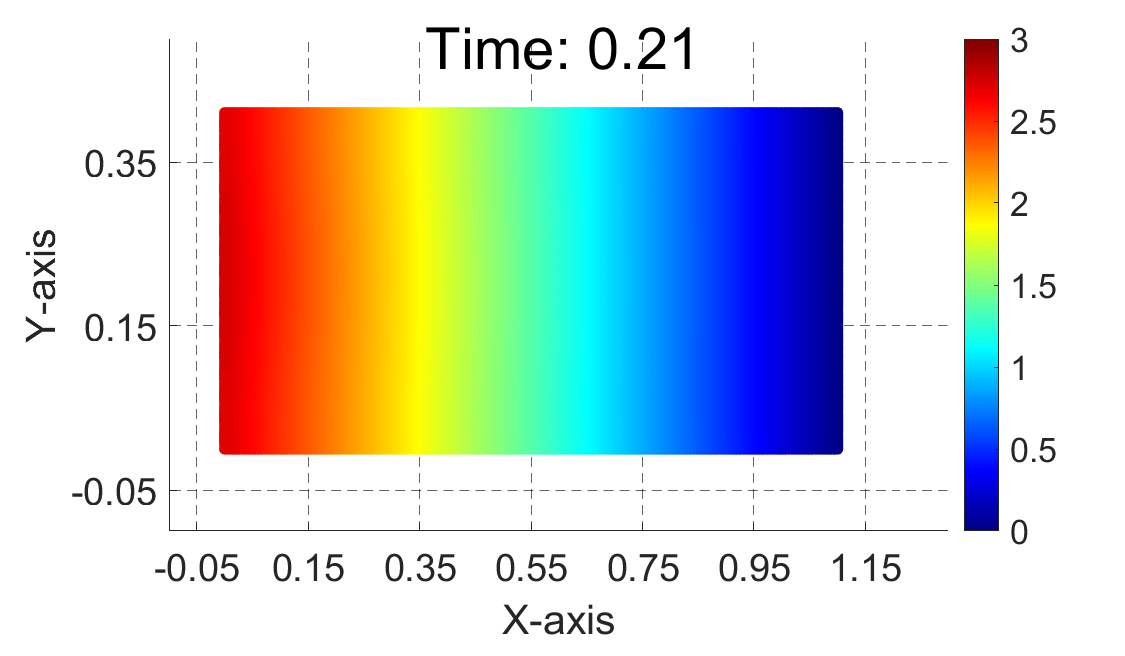}
    \includegraphics[width=0.325\textwidth]{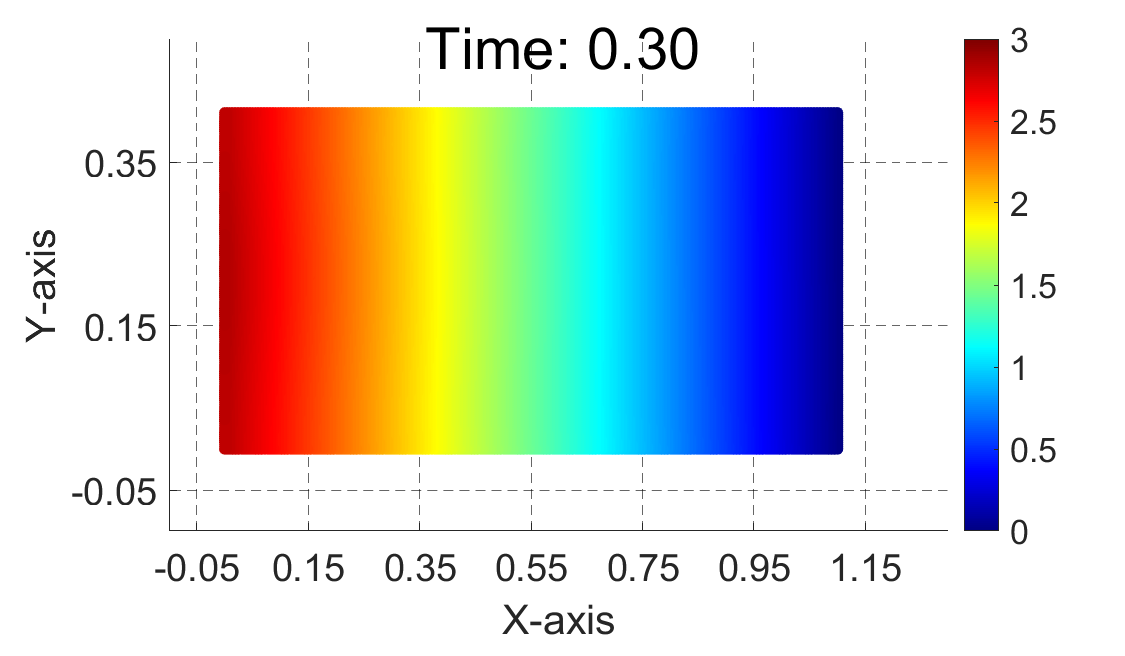}
    \includegraphics[width=0.325\textwidth]{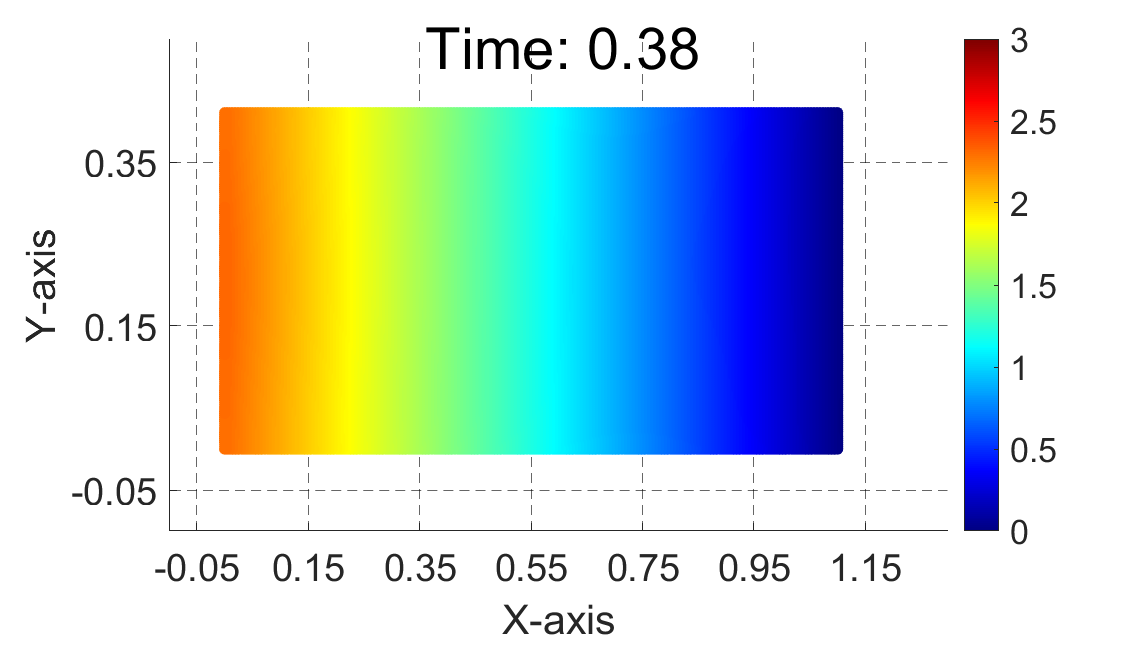}
    \includegraphics[width=0.325\textwidth]{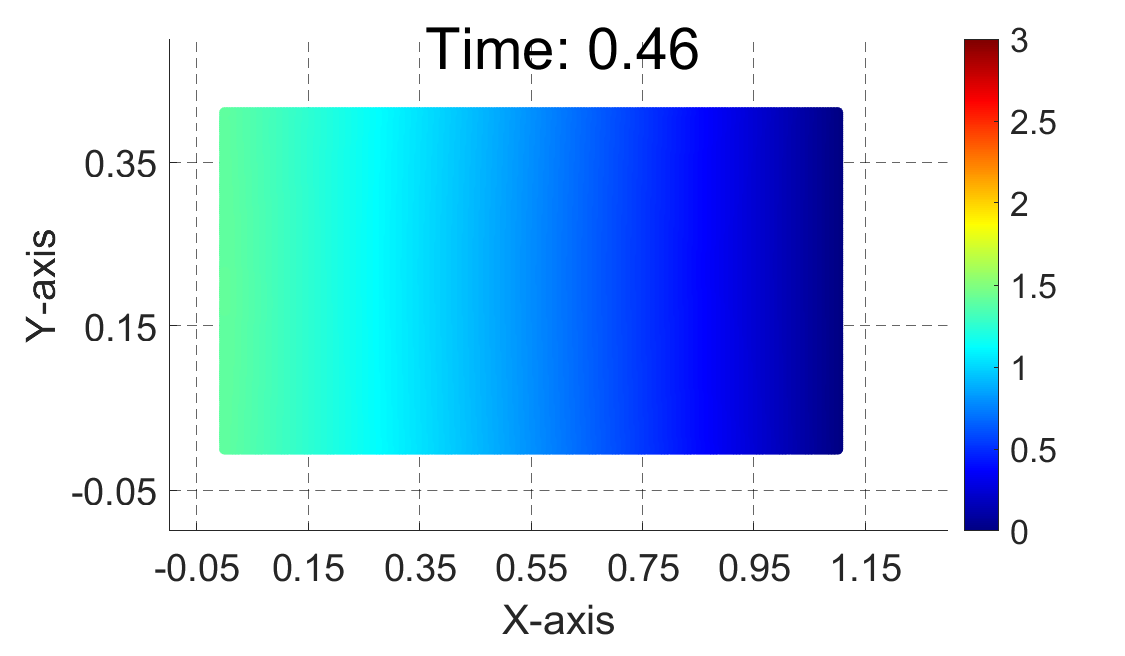}
    \caption{Predicted pressure field at various timestamps obtained from the weighted PINN model for the rectangular domain}
    \label{fig:rect_pressure_field-1}
\end{figure}

%Pressure Distribution Analysis: 
Similarly, pressure distribution plots underscored the changes in pressure gradients driving the flow, pinpointing areas of high and low pressure crucial for fluid dynamics analysis, which is depicted in Figure~\ref{fig:rect_pressure_field-1} for the rectangular domain. A pressure gradient decreasing from the inlet towards the outlet is observed, reflecting the zero-pressure boundary condition at the outlet, which drives the flow through the domain. The pressure plots exhibit pulsatile behavior, correlating with the sinusoidal inlet velocity. A cyclic increase and decrease in pressure is observed, mimicking the systolic and diastolic phases of blood flow. Comparable findings were obtained in \cite{Quarteroni_CVS00, Sarah_IJNMBE23} through conventional numerical techniques.

%Predicted Pressure Plots for rectangular domain:

The weighted PINN model performance metrics are presented below in the case of the rectangular domain. The study explored the influence of different \( \beta \) values (1, 2, 5, 10) on model accuracy and computational efficiency. The metrics of interest were final loss (accuracy indicator) and computation time, alongside the number of iterations required for convergence, which are shown in Table~\ref{table:PINN_beta-1}.

\begin{table}[H]
\centering
\caption{Summary of Performance Metrics for different \( \beta \) values in rectangular domain}
\label{tab:simplified_beta_values_rectangular}
\begin{tabular}{cccc}
\toprule
\( \beta \) & Final Loss & Comp. Time (s) & \# Iter. (Total) \\
\midrule
1 & 0.0001571 & 51338.70 & 26438 \\
2 & 0.0002181 & 45463.86 & 24081 \\
5 & 0.0002799 & 53185.37 & 27298 \\
10 & 0.0002402 & 55817.82 & 28230 \\
\bottomrule
\end{tabular}
\label{table:PINN_beta-1}
\end{table}
 
\( \beta=1 \) yielded the most accuracy while \( \beta=2 \) took the least computation time after 26,438 iterations and 24,081 iterations, respectively, optimizing the balance between accuracy and computational effort. 
 
\begin{figure}[H]
\centering
\includegraphics[width=0.45\textwidth]{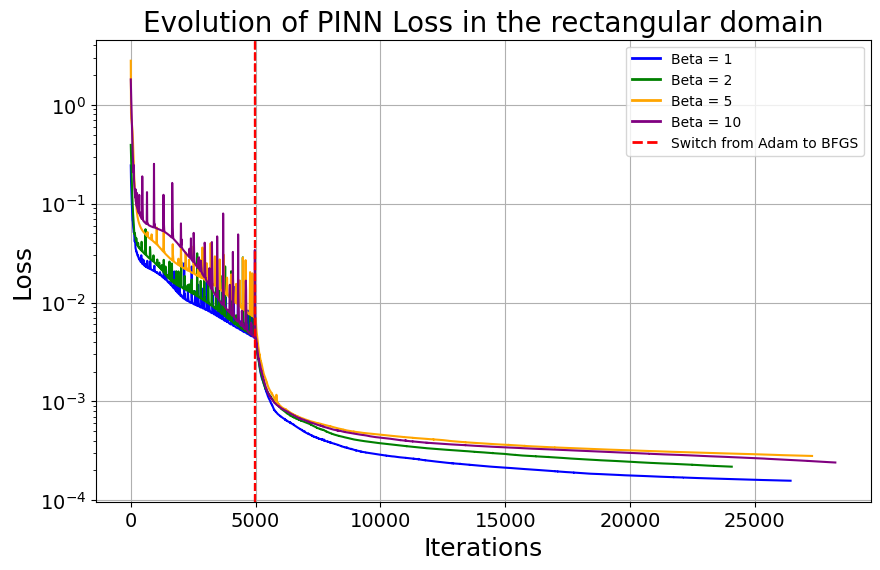}
\caption{Evolution of the loss function in the rectangular domain}
\label{fig:loss_evolution_stflow}
\end{figure}

The progression of the loss function highlighted the model’s convergence behavior across different \( \beta \) settings is depicted in Figure~\ref{fig:loss_evolution_stflow}, offering insights into performance optimization.

%\subsubsection{Semi-Circular Domain}
The weighted PINN model was also applied to a semi-circular domain to assess its performance in capturing fluid dynamics within more complex geometries.

\begin{figure}[H]
    \centering
    \includegraphics[width=0.325\textwidth]{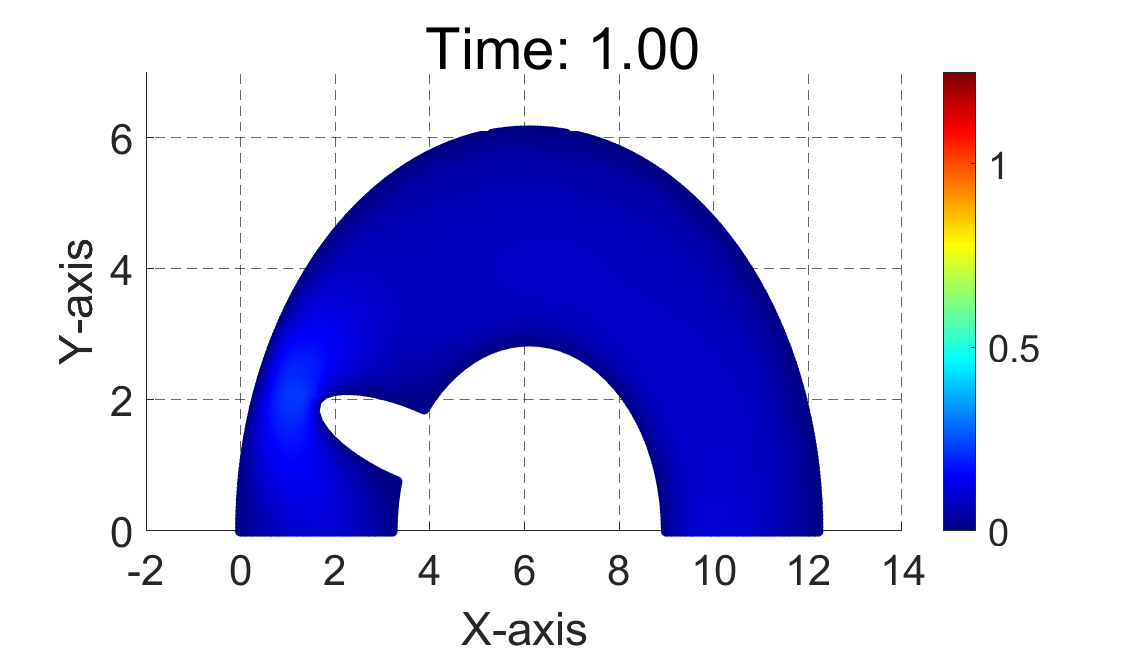}
    \includegraphics[width=0.325\textwidth]{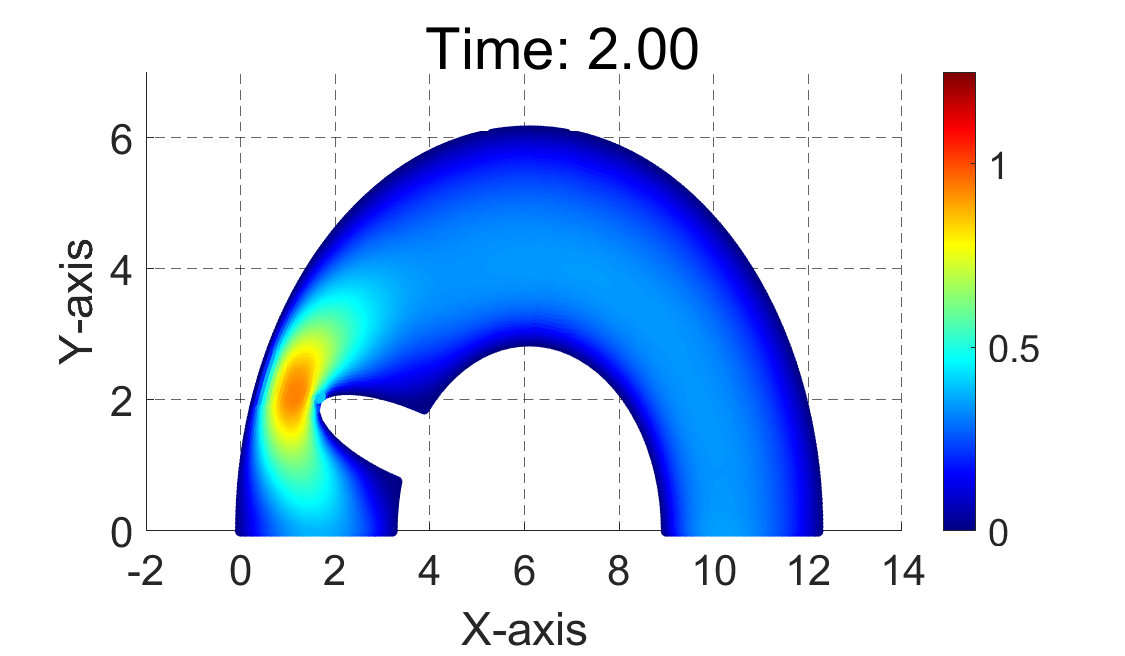}
    \includegraphics[width=0.325\textwidth]{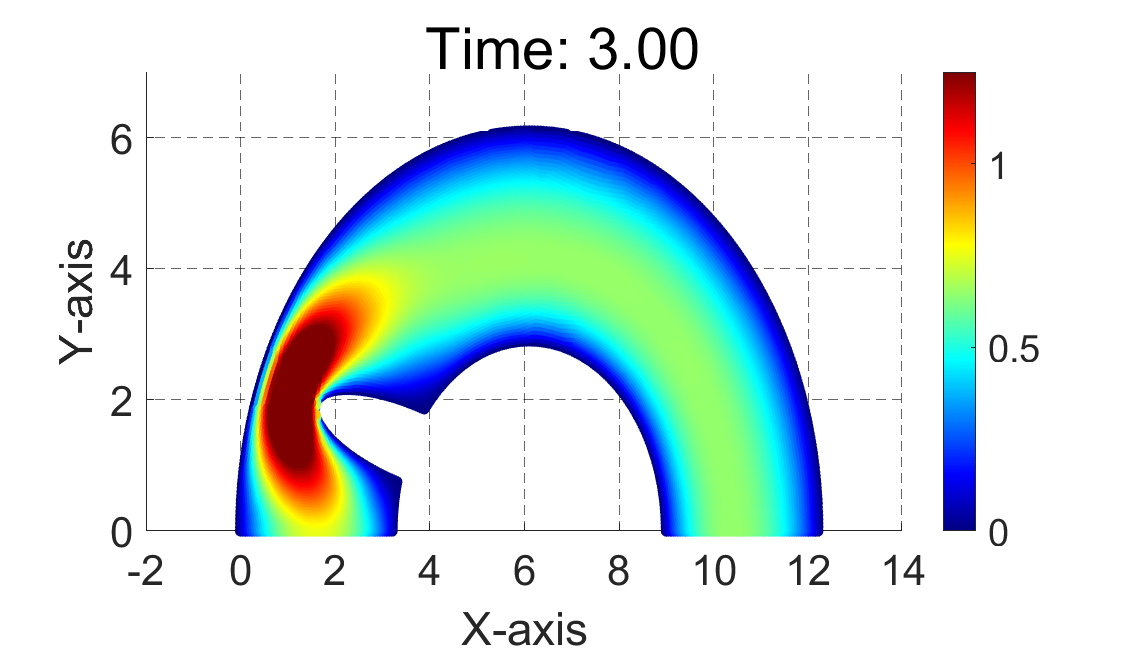}
    \includegraphics[width=0.325\textwidth]{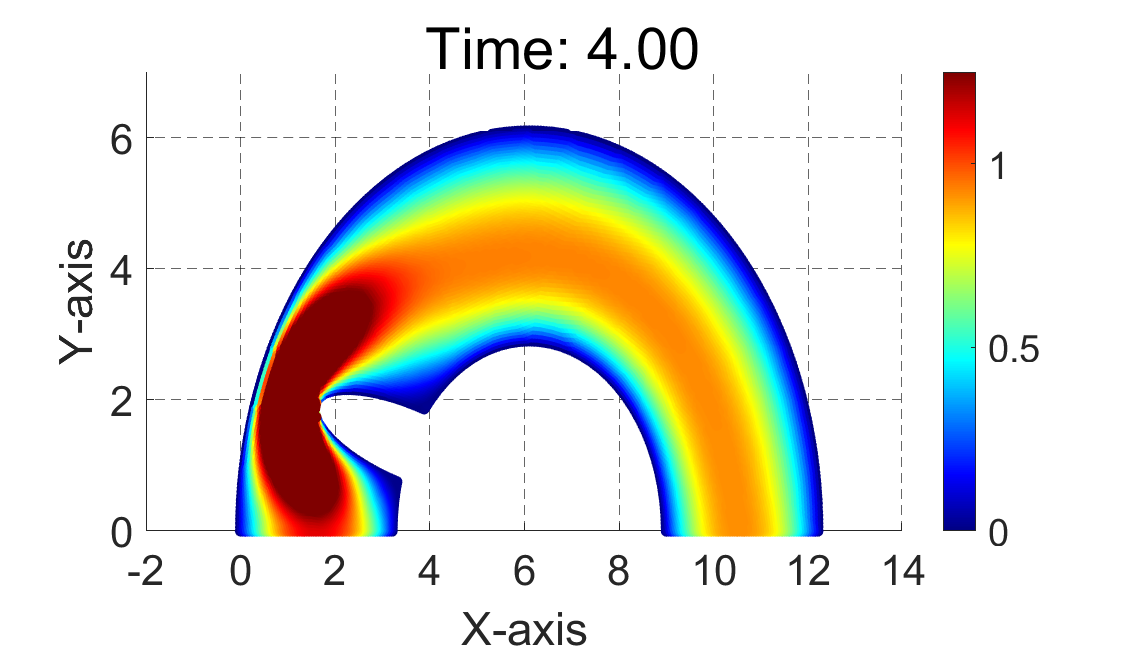}
    \includegraphics[width=0.325\textwidth]{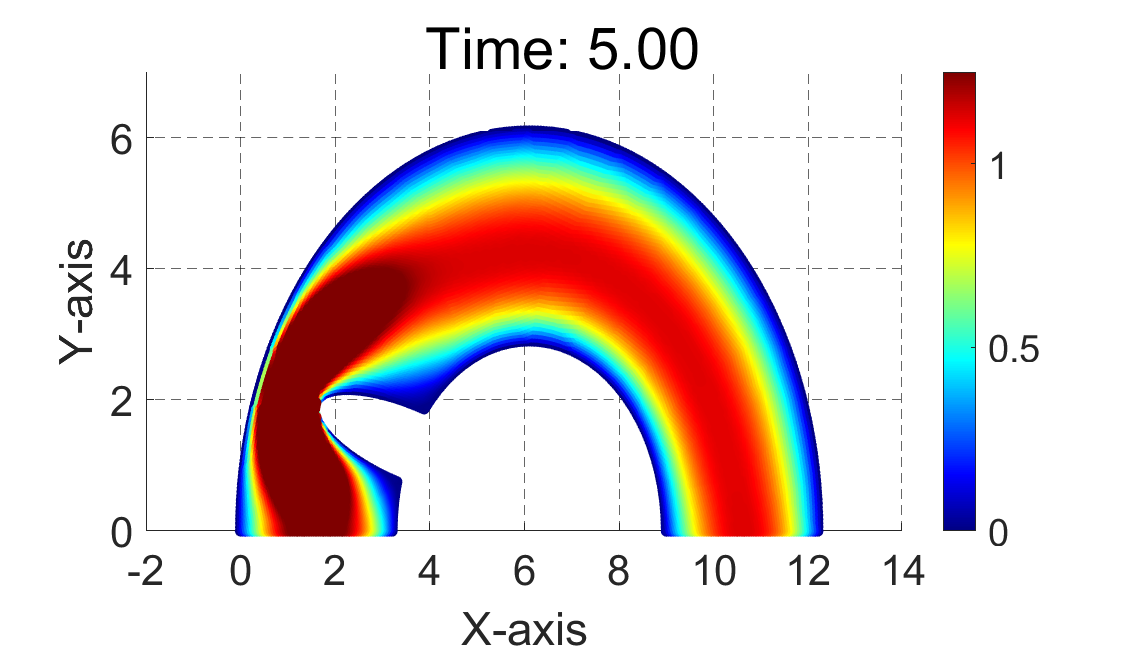}
    \includegraphics[width=0.325\textwidth]{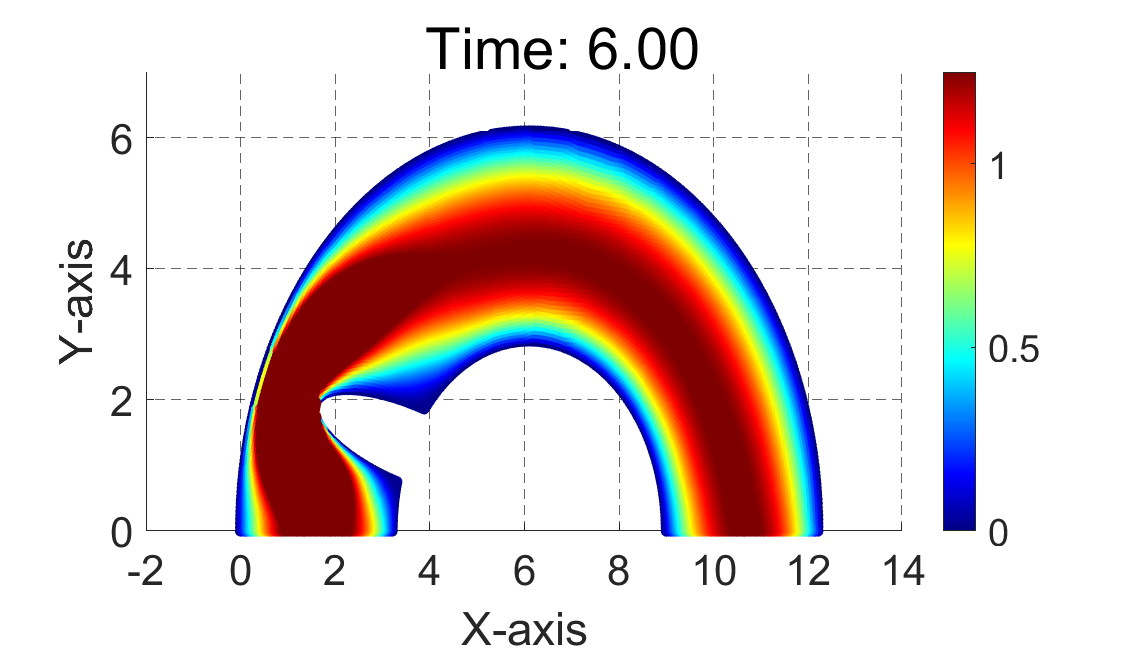}
    \caption{Predicted Velocity field at various timestamps obtained from the weighted PINN model for the semi-circular domain}
    \label{fig:semicirc_vel_field-1}
\end{figure}

%Velocity Field Analysis: 
In CFD simulations, it's a common challenge that the presence of incoming flow at open boundaries (backflow) can lead to unphysical oscillations and instability problems, even at moderate Reynolds numbers, see \cite{Cristobal_JCP16, Moghadam_CM11}. This issue arises from the convective energy entering the domain through the open boundary, which becomes problematic when the boundary velocity is unknown. Various stabilized finite element formulations have been proposed to address this problem in solving the incompressible Navier–Stokes equations. These formulations introduce stabilization terms based on the residual of a weak Stokes problem normal to the open boundary, driven by an approximate boundary pressure gradient.
In contrast, the framework of PINNs doesn't necessitate such stabilization terms to mitigate backflow instabilities. PINNs naturally handle such challenges without additional stabilization. In this study, we utilize a benchmark example problem to showcase the effectiveness of the PINNs approach. The PINNs model effectively captured the evolution of the flow field over time, highlighting significant variations in velocity across the domain, as shown in Figure~\ref{fig:semicirc_vel_field-1}. A disturbance in the flow profile is observed near the stenotic region of the artery as the velocity profile has a higher magnitude close to the origin of stenosis.

%Predicted Velocity Plots for semi-circular domain:

\begin{figure}[H]
    \centering
    \includegraphics[width=0.325\textwidth]{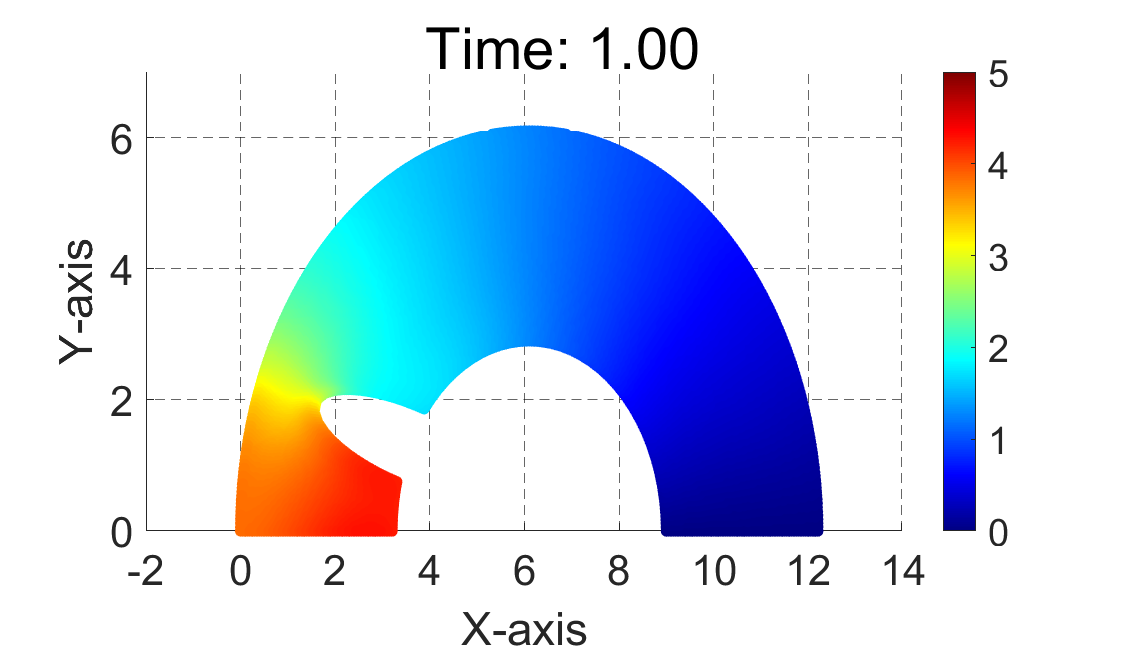}
    \includegraphics[width=0.325\textwidth]{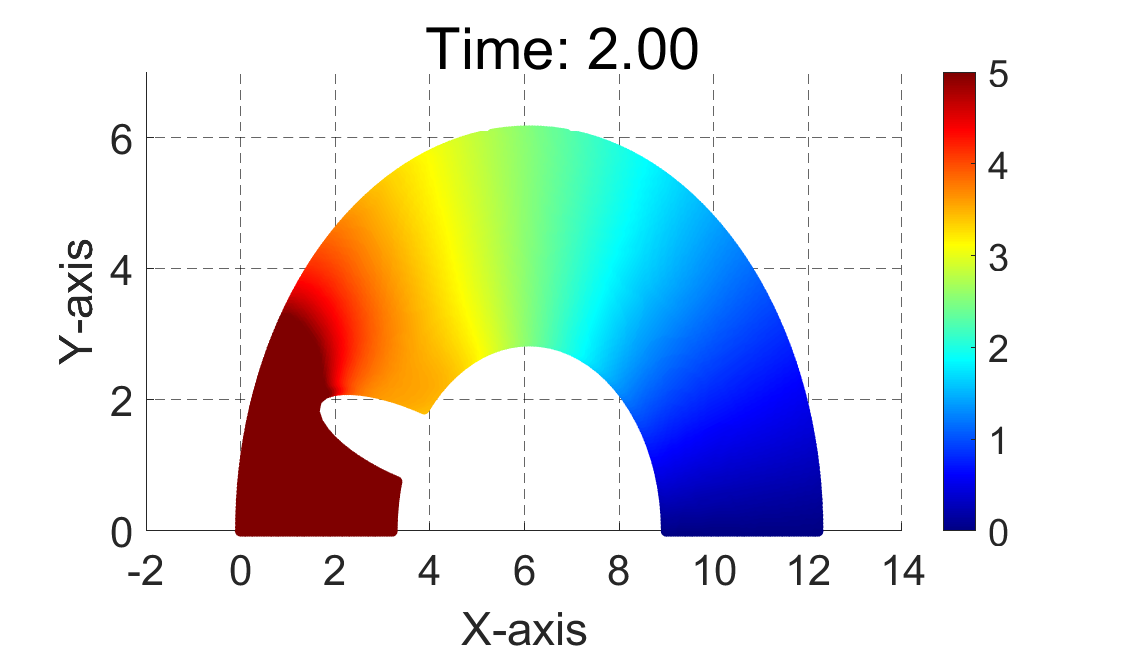}
    \includegraphics[width=0.325\textwidth]{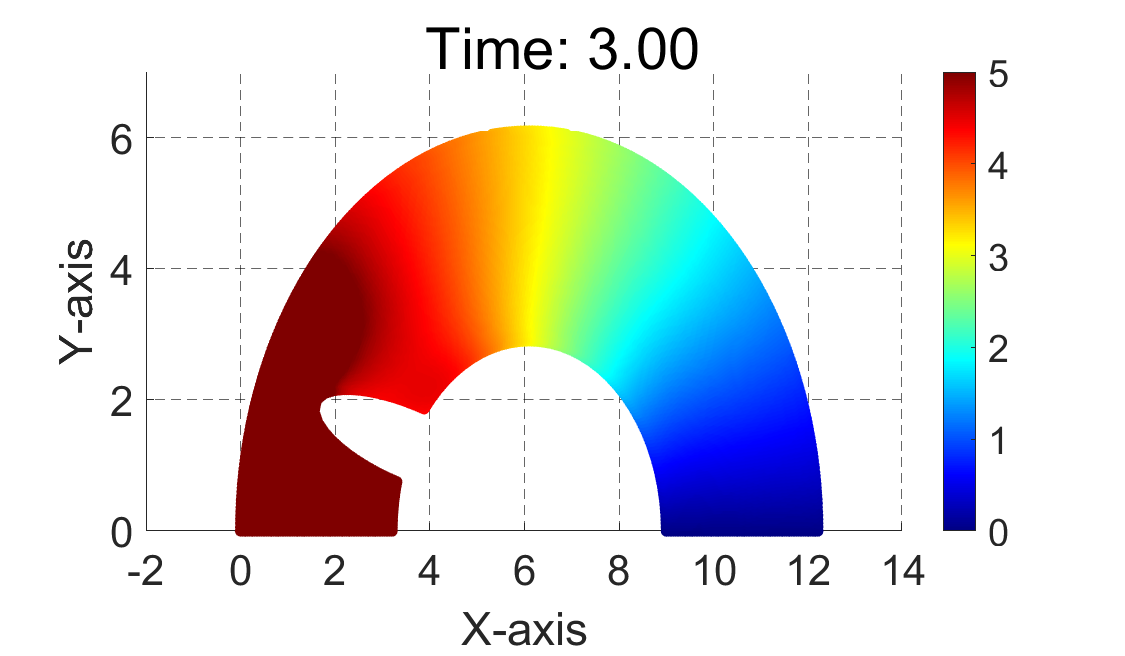}
    \includegraphics[width=0.325\textwidth]{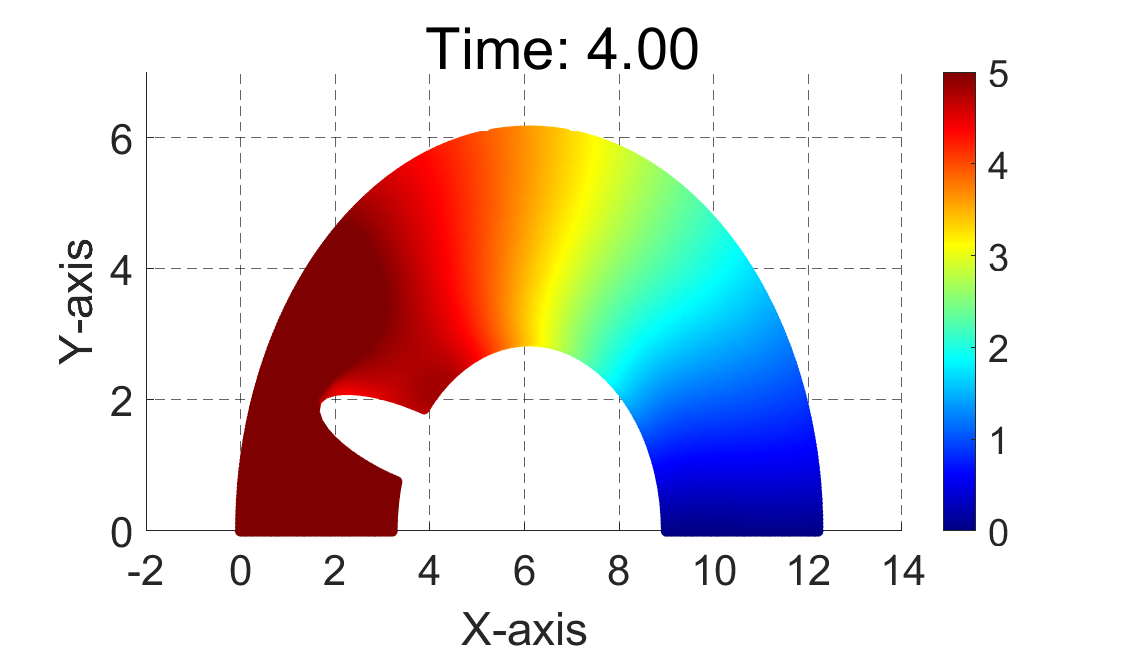}
    \includegraphics[width=0.325\textwidth]{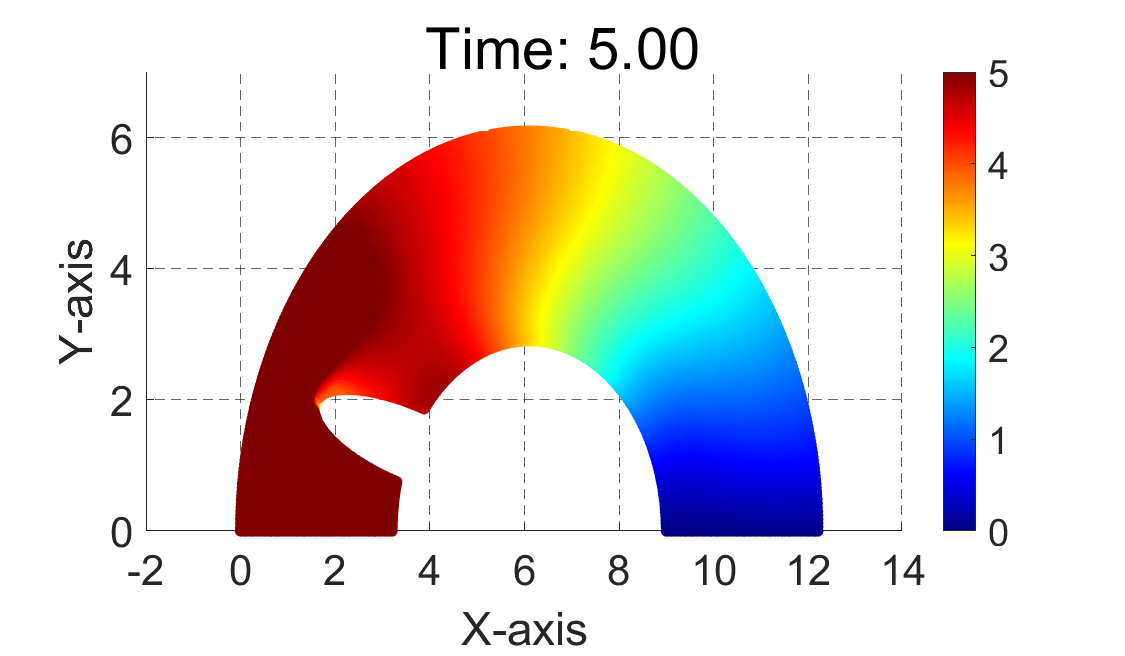}
    \includegraphics[width=0.325\textwidth]{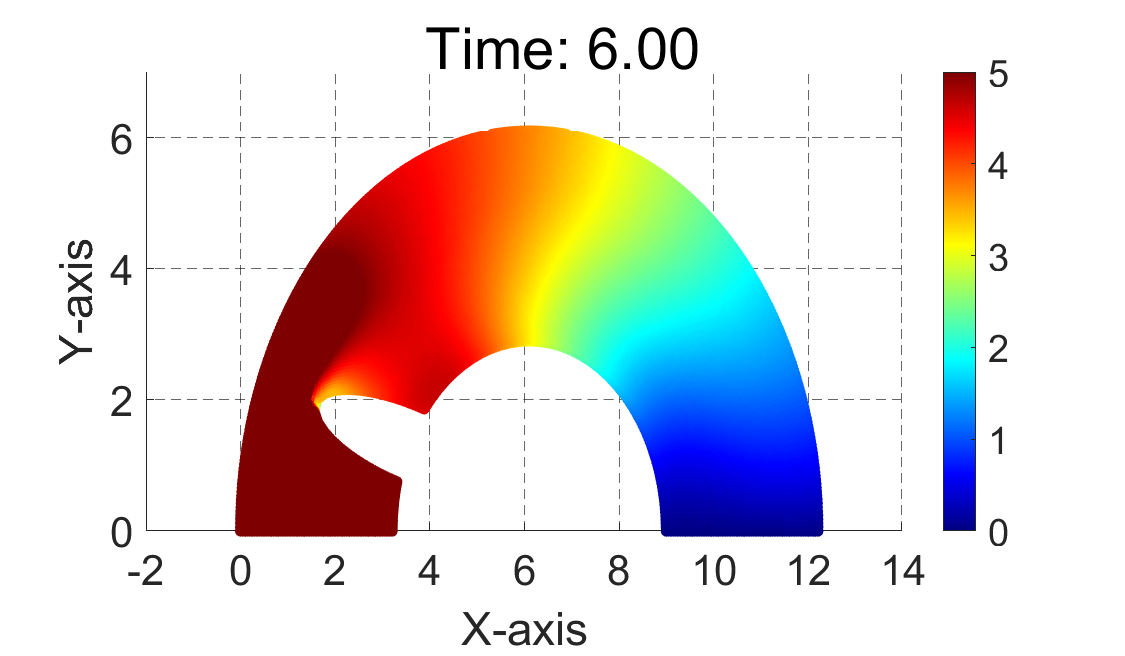}
    \caption{Predicted Pressure field at various timestamps  obtained from the weighted PINN model for the semi-circular domain}
    \label{fig:semicirc_pressure_field-1}
\end{figure}

%Pressure Distribution Analysis: 
Similarly, pressure distribution plots underscored the changes in pressure gradients driving the flow, pinpointing areas of high and low pressure crucial for fluid dynamics analysis, which is presented in Figure~\ref{fig:semicirc_pressure_field-1}. The pressure distribution also shows significant gradients around the stenotic region of the artery, with higher pressure upstream and lower downstream of the constriction. This highlights the impact of the narrowing on flow dynamics and the forces exerted on arterial walls.
Here, it's evident that we don't encounter any backflow instabilities, which stands out as the major advantage of using PINNs for addressing complex problems like this.

%Predicted Pressure Plots for semi-circular domain:

The weighted PINN model performance metrics are presented below in the case of the semi-circular domain. 
The study explored the influence of different \( \beta \) values (1, 2, 5, 10) on model accuracy and computational efficiency in Table~\ref{table:PINN_semicirc-beta-1}. 
\begin{table}[H]
\centering
\caption{Summary of Performance Metrics for different \( \beta \) values in semi-circular domain}
\label{tab:combined_beta_values_semicircular}
\begin{tabular}{cccccc}
\toprule
\( \beta \) & Final Loss & Comp. Time (s) & \# Iter. (Total) \\
\midrule
1 & 0.0003578 & 128531.72 & 16323 \\
2 & 0.0009546 & 132422.38 & 17109 \\
5 & 0.0002069 & 100092.53 & 13976 \\
10 & 0.0002037 & 123628.59 & 17389 \\
\bottomrule
\end{tabular}
\label{table:PINN_semicirc-beta-1}
\end{table}

A \( \beta \) of 10 achieved the lowest final loss after 17,389 iterations, indicating high model accuracy. Conversely, \( \beta \) of 5 was the most efficient in terms of computation time and convergence speed, with 13,976 total iterations. The progression of the loss function highlighted the model’s convergence behavior across different \( \beta \) settings, offering insights into performance optimization.
\begin{figure}[H]
\centering
\includegraphics[width=0.45\textwidth]{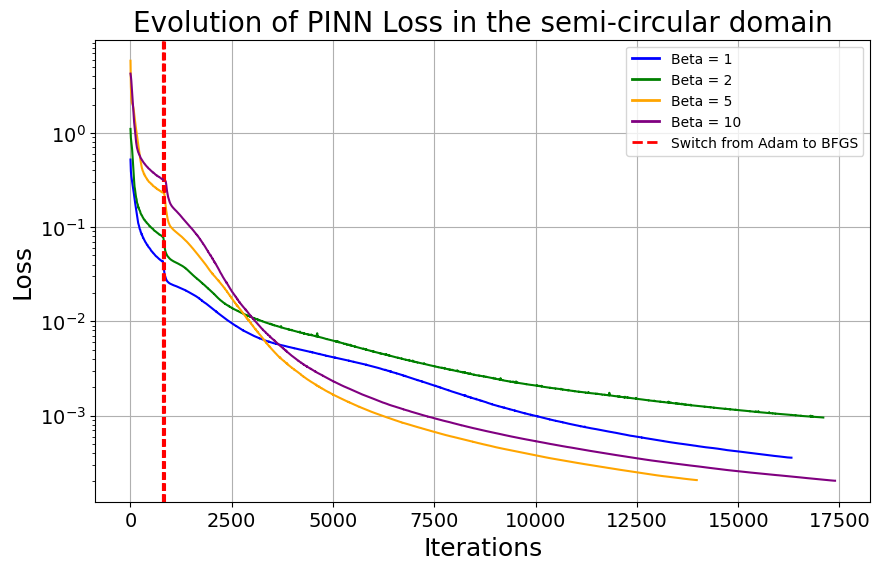}
\caption{Evolution of the loss function in the semi-circular domain}
\label{fig:loss_semicirc}
\end{figure}

%\subsection{Weighted XPINN Model}
The performance of the weighted Extended Physics-Informed Neural Network (WXPINN) model was assessed by varying the balancing coefficient \(\gamma\) and the number of subdomains \(M\) while testing with \(\beta=1\). The choice of \(\beta=1\) was taken, focusing on the model’s accuracy and computational efficiency, as that gave the best convergence in the weighted PINN case.

The study considered \(\gamma\) values of 1, 2, 5, and 10, alongside subdomain counts of 2, 3, and 4. Metrics of interest included final loss, computation time, and iterations required for convergence, providing a comprehensive view of the model's efficiency and accuracy.

Below, we demonstrate the WXPINNs efficiency w.r.t the number of subdomains:. 

\begin{itemize}
    \item Case I (\(M=2\)): In this scenario, \(\gamma=10\) recorded the shortest computation time, while \(\gamma=5\) achieved the lowest final loss, indicating the highest accuracy. See Table~\ref{table:wxpinns-rect-1} for more information.
    \item Case II (\(M=3\)): In this scenario, \(\gamma=1\) required the least computation time, and \(\gamma=5\) presented the lowest final loss, showcasing optimal accuracy. See Table~\ref{table:wxpinns-rect-1} for more information.
    \item Case III (\(M=4\)): In this scenario, \(\gamma=1\) was most efficient in terms of computation time, whereas \(\gamma=5\) maintained the lowest final loss, demonstrating superior accuracy. See Table~\ref{table:wxpinns-rect-1} for more information.
\end{itemize}

\begin{table}[H]
\centering
\caption{Summary of Performance Metrics for varying \(\gamma\) values and subdomain counts (\(M\)) in the rectangular domain for \(\beta=1\)}
\begin{tabular}{cccccc}
\toprule
\(M\) & \(\beta\) & \( \gamma \) & Final Loss & Comp. Time (s) & \# Iter. (Total) \\
\midrule
2 & 1 & 1 & 0.0001268 & 27311.14 & 20582 \\
2 & 1 & 2 & 0.0001194 & 27145.35 & 20631 \\
2 & 1 & 5 & 0.0000832 & 30256.92 & 22349 \\
2 & 1 & 10 & 0.0000869 & 26877.42 & 20733 \\
\midrule
3 & 1 & 1 & 0.0001781 & 27957.15 & 21558 \\
3 & 1 & 2 & 0.0001441 & 30425.47 & 23559 \\
3 & 1 & 5 & 0.0001276 & 28043.87 & 21337 \\
3 & 1 & 10 & 0.0001552 & 28339.32 & 21612 \\
\midrule
4 & 1 & 1 & 0.0002964 & 28353.51 & 20171 \\
4 & 1 & 2 & 0.0002528 & 31630.81 & 23064 \\
4 & 1 & 5 & 0.0001902 & 29847.30 & 20789 \\
4 & 1 & 10 & 0.0002406 & 29462.45 & 20328 \\
\bottomrule
\end{tabular}
\label{table:wxpinns-rect-1}
\end{table}

The evolution of the loss function across different configurations of the WXPINN model, varying by the number of subdomains (\(M\)) and balancing coefficient (\(\gamma\)), is visualized below. 

\begin{figure}[H]
    \centering
    \subfloat[\(M=2\) subdomains]{
        \includegraphics[width=0.45\textwidth]{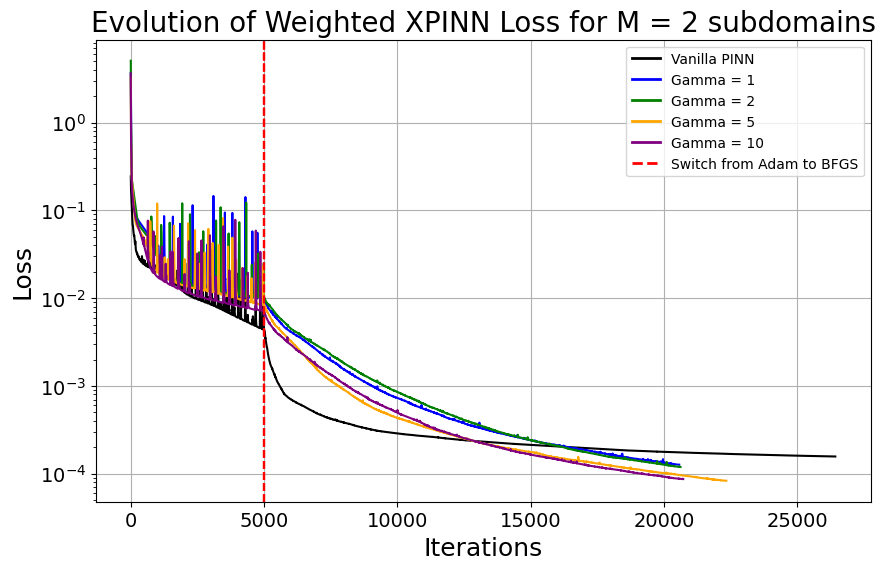}
        \label{fig:loss_m2}
    }
    \subfloat[\(M=3\) subdomains]{
        \includegraphics[width=0.45\textwidth]{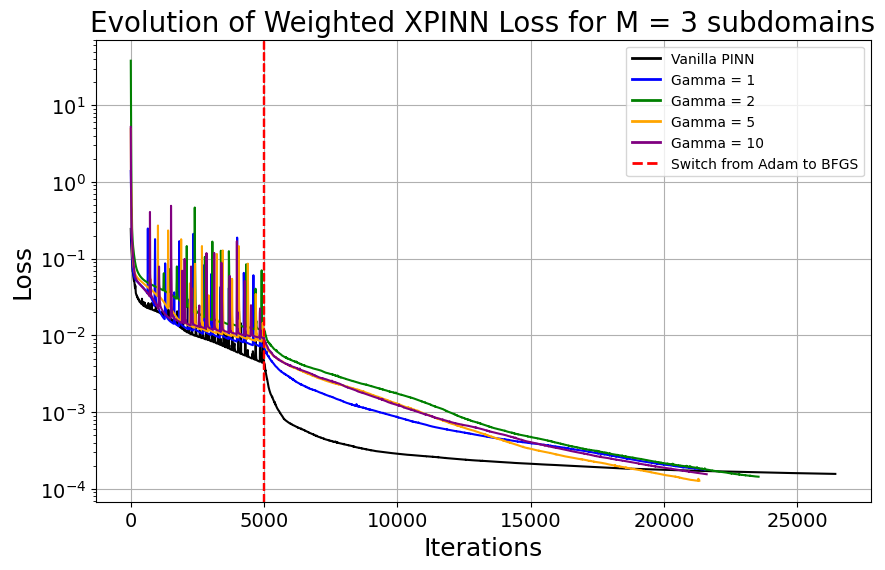}
        \label{fig:loss_m3}
    }\\
    \subfloat[\(M=4\) subdomains]{
        \includegraphics[width=0.45\textwidth]{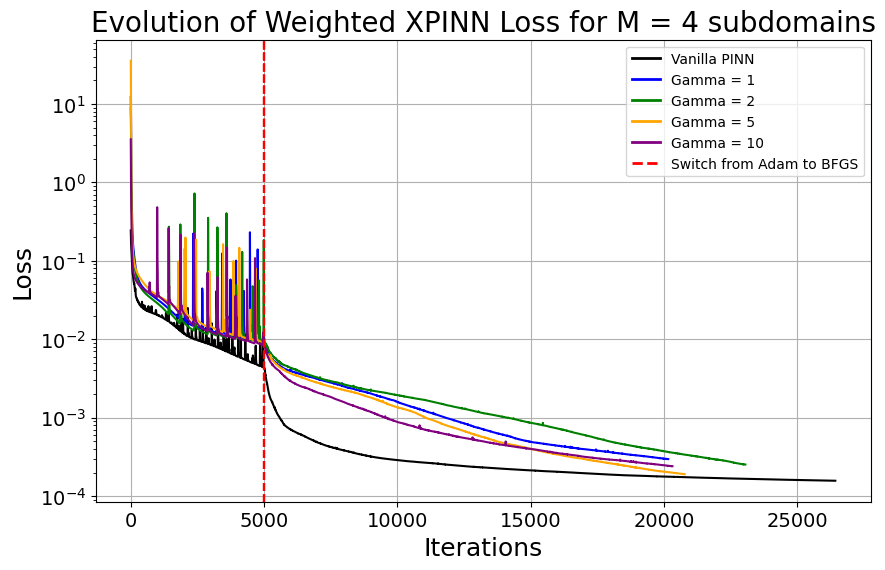}
        \label{fig:loss_m4}
    }
    \caption{Loss function evolution for the WXPINN model across different numbers of subdomains (\(M=2,3,4\)) and varying \(\gamma\) values.}
    \label{fig:consolidated_loss_evolution_xpinn}
\end{figure}

Consequently, a comparative analysis of computation times and loss evolution across various subdomain counts (\(M=2, 3, 4\)) elucidated the impact of domain partitioning on model efficiency, facilitating an understanding of accuracy-computational efficiency trade-offs.
\begin{figure}[H]
\centering
\subfloat[Computation Time]{
        \includegraphics[width=0.475\textwidth]{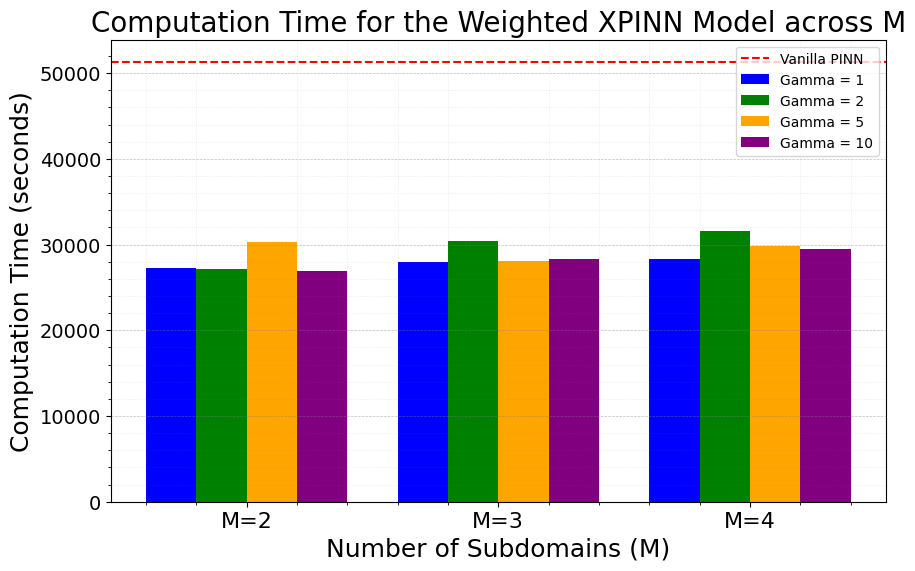}
        \label{fig:comp_time_xpinn}
    }
    \subfloat[Loss Evolution]{
        \includegraphics[width=0.45\textwidth]{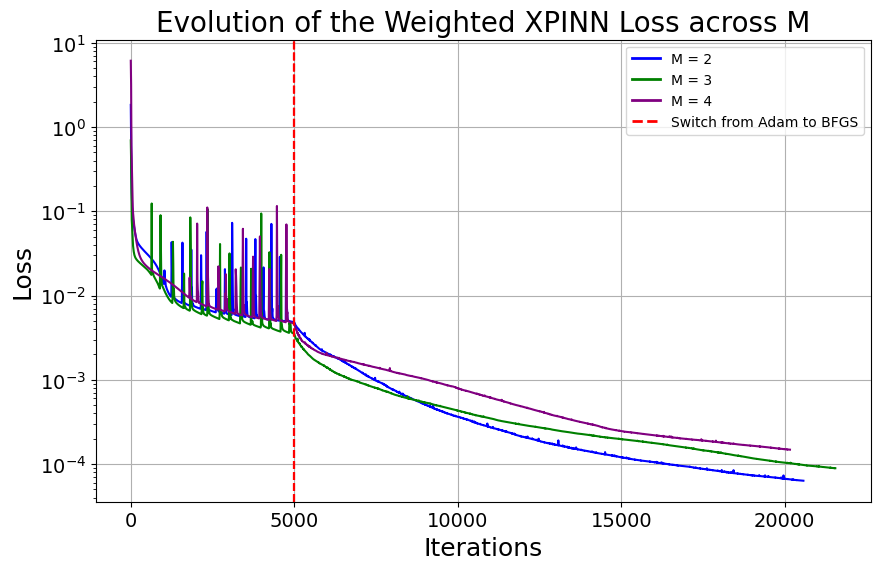}
        \label{fig:loss_m_xpinn}
    }
\caption{Comparison of Computation Time \& Loss Evolution Across Different Numbers of Subdomains (\(M=2,3,4\)) for the WXPINN model}
\label{fig:computation_time_across_subdomains_xpinn}
\end{figure}

%\subsection{Weighted CPINN Model}
The performance of the weighted Conservative Physics-Informed Neural Network (WCPINN) model was evaluated by varying the balancing coefficient \(\delta\) and the number of subdomains \(M\) while maintaining a constant \(\beta=1\) with \(\gamma\) varied as 1 and 5. The choice of \(\gamma=5\) was taken, focusing on the model’s accuracy and computational efficiency, as that gave the best convergence in the WXPINN case.

The study considered \(\delta\) values of 1, 2, 5, and 10, alongside subdomain counts of 2, 3, and 4. Metrics of interest included final loss, computation time, and iterations required for convergence, providing a comprehensive view of the model's efficiency and accuracy.

Below, we demonstrate the WCPINNs efficiency w.r.t the number of subdomains:

\begin{itemize}
    \item Case I (\(M=2\)): In this scenario, for \(\gamma=1\), \(\delta=10\)
     recorded the shortest computation time, while \(\delta=2\) achieved the lowest final loss, and for \(\gamma=5\), \(\delta=5\) recorded the shortest computation time, while \(\delta=2\) achieved the lowest final loss, indicating the highest accuracy. See Table~\ref{table:wcpinns-rect-1} for more information.
    \item Case II (\(M=3\)): In this scenario, for \(\gamma=1\), \(\delta=5\) required the least computation time, and \(\delta=2\) presented the lowest final loss, and for \(\gamma=5\), \(\delta=1\) required the least computation time, and \(\delta=2\) presented the lowest final loss, showcasing optimal accuracy. See Table~\ref{table:wcpinns-rect-1} for more information.
    \item Case III (\(M=4\)): In this scenario, for \(\gamma=1\), \(\delta=10\) was most efficient in terms of computation time, whereas \(\delta=2\) maintained the lowest final loss, and for \(\gamma=5\), \(\delta=5\) was most efficient in terms of computation time, whereas \(\delta=2\) maintained the lowest final loss, demonstrating superior accuracy. See Table~\ref{table:wcpinns-rect-1} for more information.
\end{itemize}

\begin{table}[H]
\centering
\caption{Summary of Performance Metrics for varying \(\delta\) values and subdomain counts (\(M\)) in the rectangular domain for \(\gamma=1\) and \(\gamma=5\)}
\begin{tabular}{ccccccccc} 
\toprule
\multirow{2}{*}{\(M\)} & \multirow{2}{*}{\(\delta\)} & \multicolumn{2}{c}{Final Loss} & \multicolumn{2}{c}{Comp. Time (s)} & \multicolumn{2}{c}{\# Iter. (Total)} \\
\cmidrule(lr){3-4} \cmidrule(lr){5-6} \cmidrule(lr){7-8}
& & \(\gamma=1\) & \(\gamma=5\) & \(\gamma=1\) & \(\gamma=5\) & \(\gamma=1\) & \(\gamma=5\) \\
\midrule
2 & 1 & 0.0000834 & 0.0001167 & 26825.72 & 27398.77 & 20006 & 20444 \\
2 & 2 & 0.0000750 & 0.0000693 & 30754.57 & 28387.97 & 23268 & 21068 \\
2 & 5 & 0.0000965 & 0.0000723 & 29168.68 & 26062.06 & 22780 & 20240 \\
2 & 10 & 0.0001355 & 0.0001567 & 26260.77 & 30778.42 & 19424 & 23276 \\
\midrule
3 & 1 & 0.0001898 & 0.0001689 & 30810.61 & 26391.59 & 23689 & 19864 \\
3 & 2 & 0.0001107 & 0.0001056 & 28523.03 & 27621.32 & 20824 & 20549 \\
3 & 5 & 0.0002403 & 0.0001206 & 26233.19 & 28073.90 & 19561 & 21535 \\
3 & 10 & 0.0002224 & 0.0001667 & 29465.80 & 26601.74 & 20719 & 20153 \\
\midrule
4 & 1 & 0.0003783 & 0.0003722 & 31109.31 & 29089.39 & 21589 & 20204 \\
4 & 2 & 0.0002296 & 0.0002105 & 30765.06 & 29663.42 & 21036 & 19551 \\
4 & 5 & 0.0003369 & 0.0002911 & 28601.59 & 29013.70 & 19953 & 19915 \\
4 & 10 & 0.0003777 & 0.0002347 & 27822.58 & 29745.53 & 19097 & 19565 \\
\bottomrule
\end{tabular}
\label{table:wcpinns-rect-1}
\end{table}

The evolution of the loss function across different configurations of the WCPINN model, varying by the number of subdomains (\(M\)) and balancing coefficient (\(\delta\)), is visualized below. 

\begin{figure}[H] 
    \centering
    \subfloat[\(M=2\) subdomains (\(\gamma=1\))]{
        \includegraphics[width=0.45\textwidth]{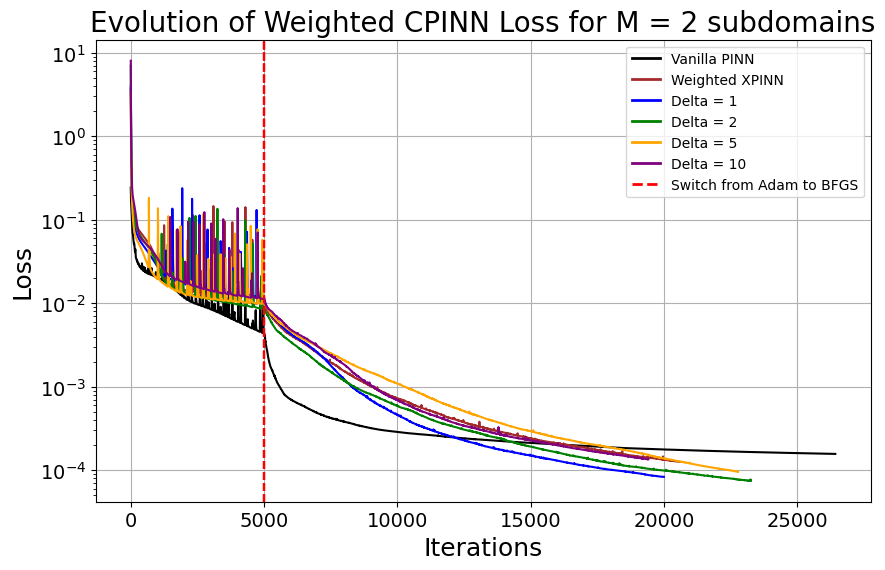}
        \label{fig:loss_m2c1}
    }
    \subfloat[\(M=2\) subdomains (\(\gamma=5\))]{
        \includegraphics[width=0.45\textwidth]{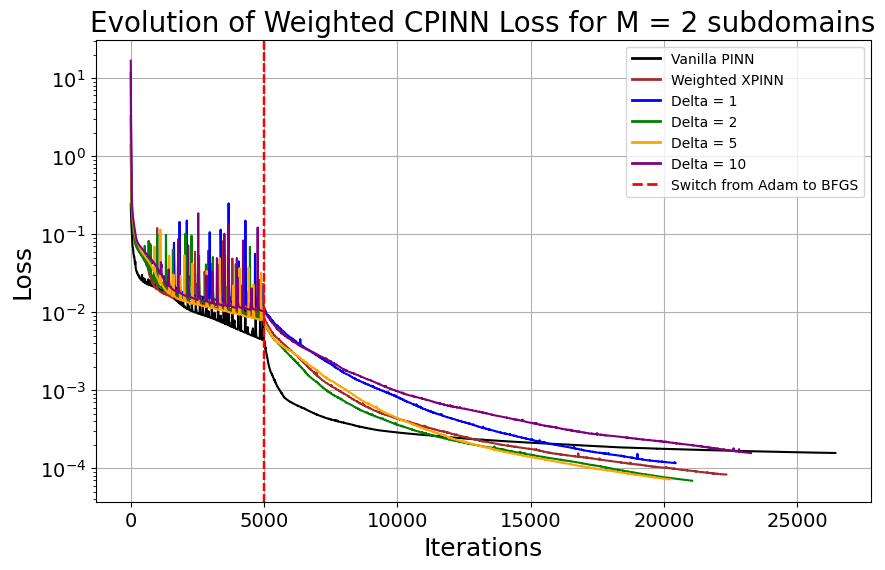}
        \label{fig:loss_m2c5}
    }\\
    \subfloat[\(M=3\) subdomains (\(\gamma=1\))]{
        \includegraphics[width=0.45\textwidth]{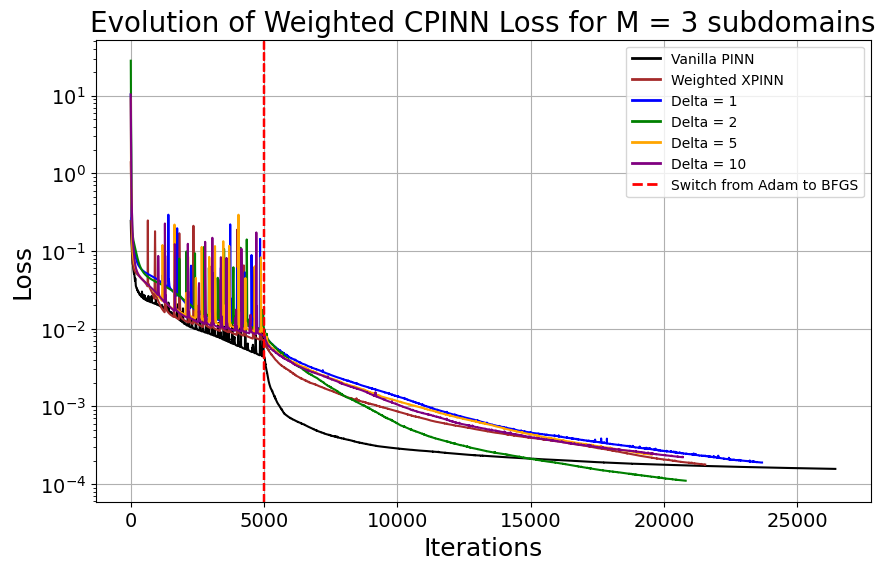}
        \label{fig:loss_m3c1}
    }
    \subfloat[\(M=3\) subdomains (\(\gamma=5\))]{
        \includegraphics[width=0.45\textwidth]{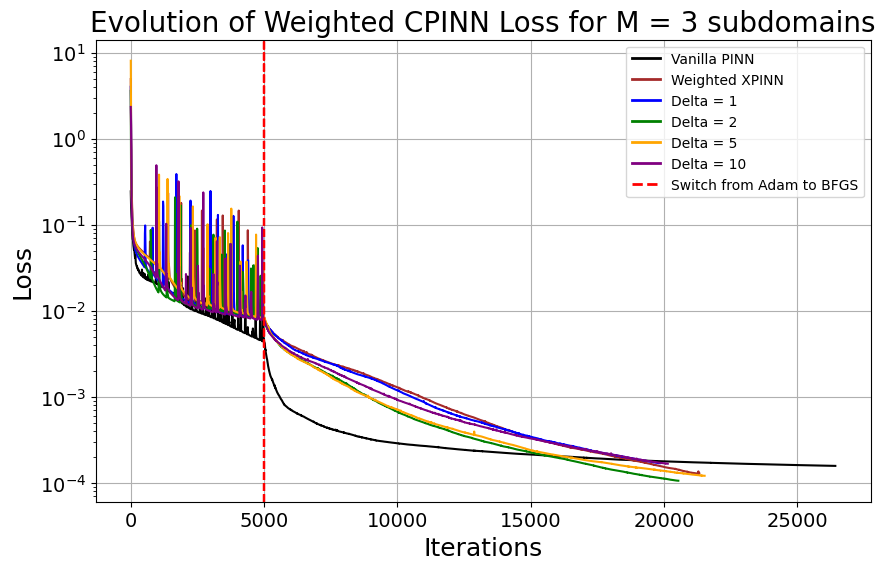}
        \label{fig:loss_m3c5}
    }
    \caption{Loss function evolution for the WCPINN model, part 1.}
    \label{fig:consolidated_loss_evolution_cpinn_part1}
\end{figure}

\begin{figure}[H] 
    \ContinuedFloat 
    \centering
    \subfloat[\(M=4\) subdomains (\(\gamma=1\))]{
        \includegraphics[width=0.45\textwidth]{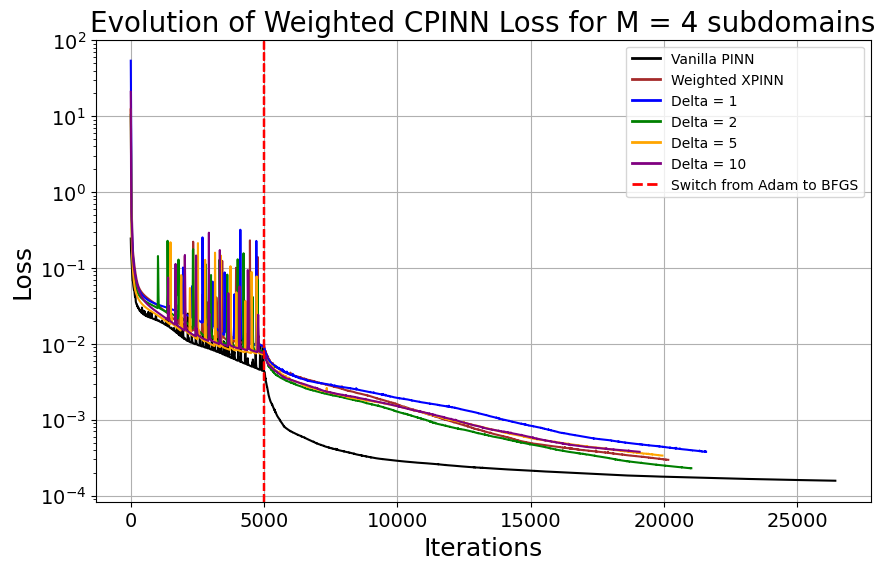}
        \label{fig:loss_m4c1}
    }
    \subfloat[\(M=4\) subdomains (\(\gamma=5\))]{
        \includegraphics[width=0.45\textwidth]{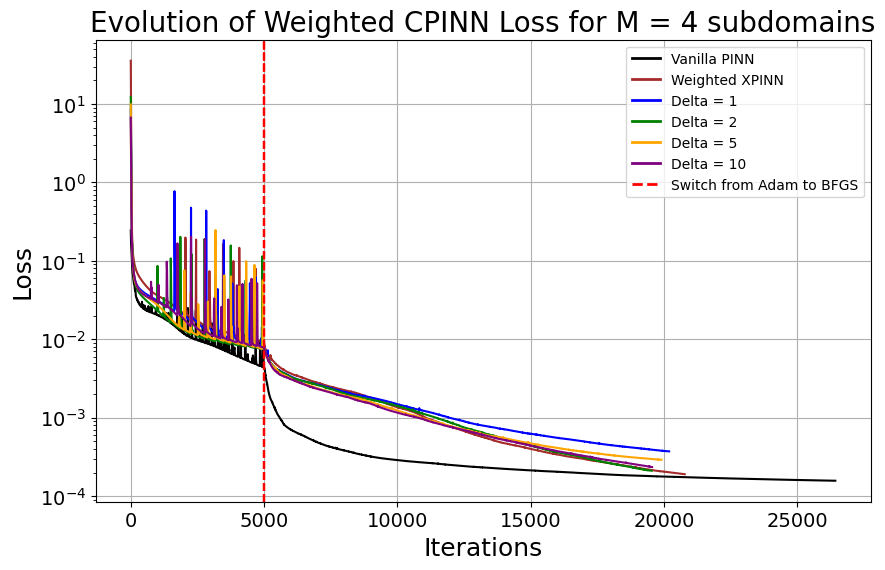}
        \label{fig:loss_m4c5}
    }
    \caption{Loss function evolution for the WCPINN model, part 2 (Continued).}
    \label{fig:consolidated_loss_evolution_cpinn_part2}
\end{figure}

A comparative evaluation of computation times further elucidates the impact of domain partitioning on computational efficiency, assisting in navigating the trade-offs between simulation accuracy and computational demand.
\begin{figure}[H]
\centering
\subfloat[Computation Time (\(\gamma=1\))]{
        \includegraphics[width=0.475\textwidth]{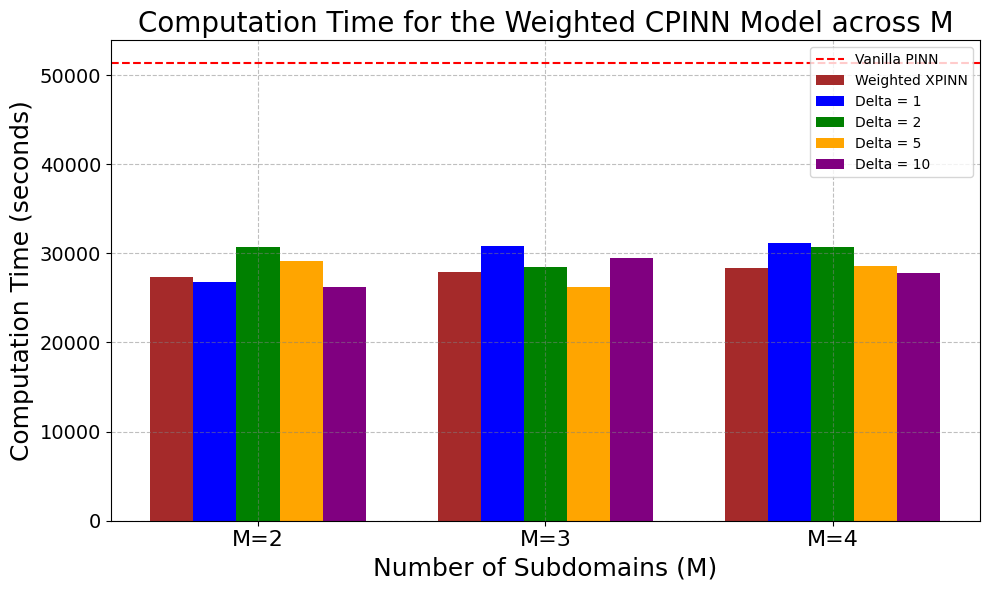}
        \label{fig:comp_time_cpinn_1}
    }
\subfloat[Computation Time (\(\gamma=5\))]{
        \includegraphics[width=0.475\textwidth]{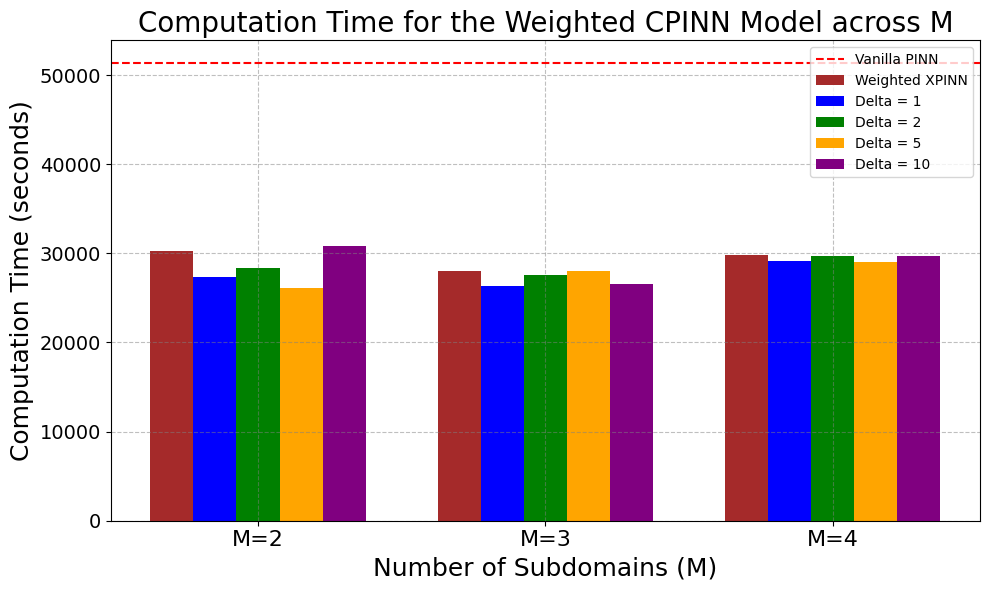}
        \label{fig:comp_time_cpinn_5}
    }   \\ 
\subfloat[Loss Evolution (\(\gamma=1\))]{
        \includegraphics[width=0.45\textwidth]{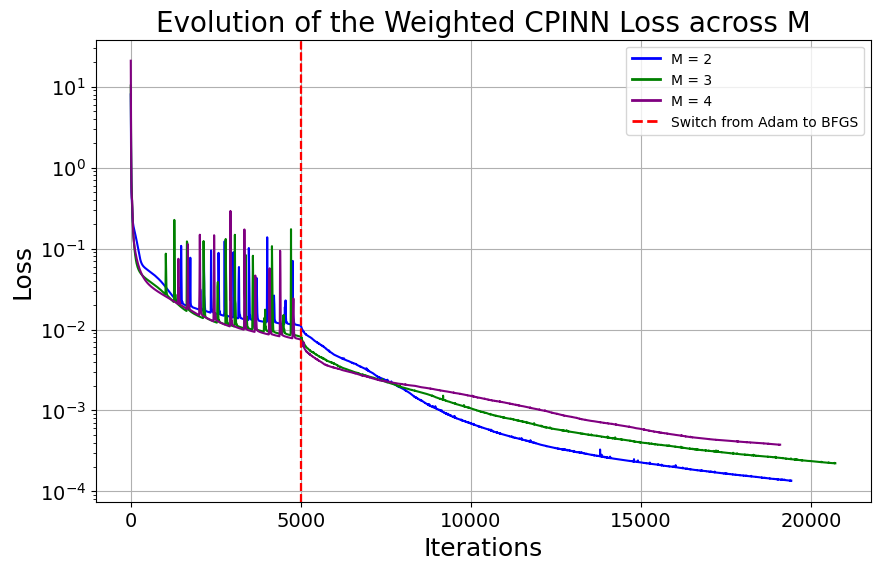}
        \label{fig:loss_m_cpinn_1}
    }
\subfloat[Loss Evolution (\(\gamma=5\))]{
        \includegraphics[width=0.45\textwidth]{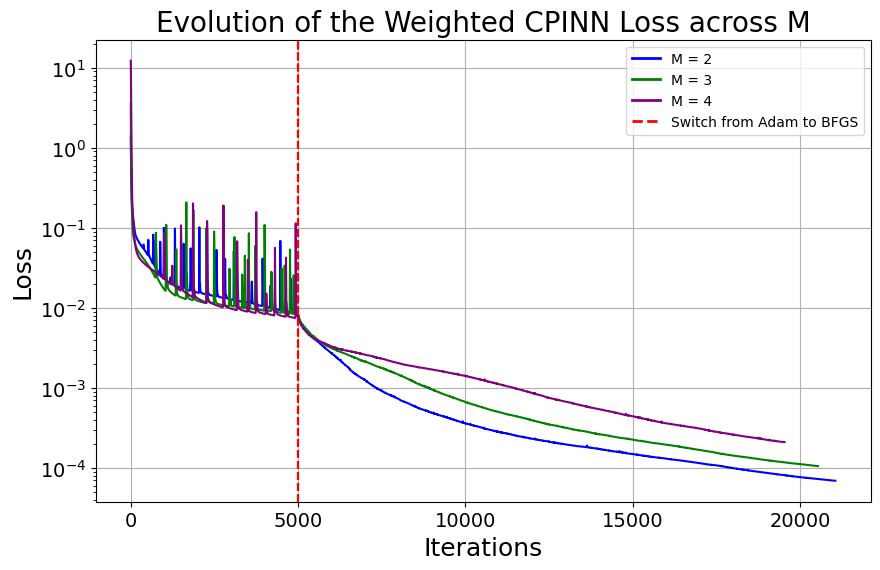}
        \label{fig:loss_m_cpinn_5}
    }
\caption{Comparison of Computation Time \& Loss Evolution Across Different Numbers of Subdomains (\(M=2,3,4\)) for the WCPINN model}
\label{fig:computation_time_across_subdomains_cpinn}
\end{figure}

\section{Conclusion}
This study rigorously explores the capabilities of Physics-Informed Neural Networks (PINNs), including their advanced variants as weighted Extended PINN (WXPINN) and weighted Conservative PINN (WCPINN), to solve the blood flow model equations within cardiovascular systems. The study demonstrates that strategic subdomain partitioning and the optimal selection of balancing coefficients (\(\beta\), \(\gamma\), \(\delta\)) are essential for enhancing computational efficiency while preserving high model accuracy. The adaptability of these models to diverse geometries underscores their potential for broader applications in simulating real-world fluid dynamics scenarios. Our PINNs simulation results clearly demonstrate the absence of backflow instabilities, highlighting a significant advantage of utilizing PINNs over traditional numerical methods for tackling complex problems like this.

However, the performance of PINN-based models heavily relies on hyperparameter configuration, necessitating extensive tuning for optimal settings. While domain decomposition improves generalization across subdomains, it introduces a tradeoff: simplifying each segment reduces required learning complexity but also decreases available training data per subdomain, potentially increasing overfitting risk and diminishing model generalizability. To address these challenges, we investigated WXPINN and WCPINN models incorporating weighted methodologies to optimize the distribution of training data and computational resources across subdomains for solving Navier-Stokes equations on complex domains. Based on presented numerical results and strategic adaptations, WXPINN and WCPINN provide a general approach to managing inherent tradeoffs in domain-decomposed neural network training, thereby enhancing the accuracy and applicability of PINN simulations in complex fluid dynamics.

Future endeavors should prioritize improving the generalization abilities of weighted PINNs across parallel computing frameworks and GPU accelerators. This enhancement aims to enhance computational efficiency, facilitating the scaling of these models on complex geometries and achieving real-time performance.

\section*{Acknowledgement} 
We gratefully acknowledge the financial support under R\&D projects leading to HPC Applications (DST/NSM/R\&D\_HPC\_Applications/Extension/2023/33), National Supercomputing Mission, India, and we greatly acknowledge the support for high-performance computing facilities at the Padmanabha cluster, IISER Thiruvananthapuram, India. Additionally, we extend our gratitude to Miss Soundarya Sarathi for her initial code implementation.

\bibliographystyle{ieeetr}
\bibliography{ref.bib}

\end{document}